\documentclass[11pt]{article}

\usepackage{amsthm,amsmath,amssymb,bm}
\usepackage{graphicx,color,epsfig}
\usepackage[margin=1in]{geometry}

\newcommand{\N}{\mathbb{N}}
\newcommand{\R}{\mathbb{R}}
\newcommand{\T}{\mathbb{T}}
\newcommand{\Z}{\mathbb{Z}}

\newtheorem{Thm}{Theorem}[section]
\newtheorem{Pro}[Thm]{Proposition}
\newtheorem{Lem}[Thm]{Lemma}
\newtheorem{Claim}[Thm]{Claim}
\newtheorem{Exa}[Thm]{Example}
\newtheorem{Rem}[Thm]{Remark}
\newtheorem{Def}[Thm]{Definition}

\begin{document}

\title{Symmetry-breaking-induced loss of ergodicity in maps of the simplex with inversion symmetry}
\author{Bastien Fernandez and Eric Vernier}
\date{}
\maketitle
\begin{center}
Laboratoire de Probabilit\'es, Statistique et Mod\'elisation\\
CNRS - Univ. Paris Cit\'e -  Sorbonne Univ.\\
Paris, France\\
fernandez@lpsm.paris and vernier@lpsm.paris
\end{center}

\begin{abstract}
Motivated by proving the loss of ergodicity in expanding systems of piecewise affine coupled maps with arbitrary number of units, all-to-all coupling and inversion symmetry, we provide ad-hoc substitutes - namely inversion-symmetric maps of the simplex with arbitrary number of vertices - that exhibit several asymmetric absolutely continuous invariant measures when their expanding rate is sufficiently small. In a preliminary study, we consider arbitrary maps of the multi-dimensional torus with permutation symmetries. Using these symmetries, we show that the existence of multiple invariant sets of such maps can be obtained from their analogues in some reduced maps of a smaller phase space. For the coupled maps, this reduction yields inversion-symmetric maps of the simplex. The subsequent analysis of these reduced maps show that their systematic dynamics is intractable because some essential features vary with the number of units; hence the substitutes which nonetheless capture the coupled maps common characteristics. The construction itself is based on a simple mechanism for the generation of asymmetric invariant union of polytopes, whose basic principles should extend to a broad range of maps with permutation and inversion symmetries.  
\end{abstract}

\leftline{\small\today.}

\section{Introduction}
\subsection{Background and motivations}
Systems of coupled maps were introduced as discrete time models for the dynamics of collective systems of interacting units \cite{K93}. They have revealed a rich phenomenology depending on the individual dynamics and on the coupling type and strength. Part of this phenomenology has been proved from a rigorous mathematical point of view \cite{CF05}. 

In particular, the most common result in the chaotic (expanding or hyperbolic) setting is the existence of a unique absolutely continuous invariant measure ({\bf acim}) when the coupling strength is sufficiently weak. Such uniqueness follows from perturbation arguments at the uncoupled limit, see e.g.\ \cite{KK92} for piecewise expanding coupled maps with a finite number of units. 
For systems with infinitely many units, similar uniqueness statements have been proved in the case of nearest-neighbour or exponentially decaying coupling, see for instance \cite{BS88,JP98,KL06,PS88}. These results have been considered as the analogue in the deterministic setting of the theory of dynamical systems, of the uniqueness of the high temperature phase in particle systems of statistical mechanics, and especially in the Ising model. 

The analogy with particle systems suggests that, when the coupling strength increases sufficiently, uniqueness of the acim (and hence ergodicity) should be lost via some analogue of a symmetry-breaking-induced phase transition \cite{BS88}. Yet, the nature of this transition and its outcomes have long been debated in the community \cite{CF05}. In particular, outside the weak coupling regime, the features of the symbolic dynamics on which the thermodynamics formalism and the subsequent theory of phase transitions are grounded, are not usually known with enough detail. In the setting of infinite lattices, exceptions have been provided by ad-hoc examples inspired by Toom's cellular automata that exhibit standard phase transitions \cite{BK06,GM00}. 

Loss of ergodicity upon sufficient increase of the coupling strength also occurs in systems with finitely many units. This is particularly the case of 
the family $\{F_{N,\epsilon}\}_{\epsilon\in [0,\frac12)}$ of maps of the $N$-dimensional torus $\T^N$ defined by \cite{KY10}\footnote{We use $[1,N]$ to denote the collection of the first $N$ natural integers.}  
\[
(F_{N,\epsilon}\mathrm{u})_i= 2\left(\mathrm{u}_i+\tfrac{\epsilon}N\sum_{j=1}^Ng(\mathrm{u}_j-\mathrm{u}_i)\right)\ \text{mod}\ 1,\ \forall i\in [1,N],\ \mathrm{u}=(\mathrm{u}_i)_{i=1}^N\in\T^N,
\]
where\footnote{The symbol $\lfloor\cdot\rfloor$ denotes the floor function.}
\[ 
g(\mathrm{u})=\left\{\begin{array}{ccl}
\mathrm{u}-\lfloor \mathrm{u}+\tfrac12\rfloor&\text{if}&\mathrm{u}\not\in \tfrac12\ \text{mod}\ 1\\
\mathrm{u}&\text{if}&\mathrm{u}\in \tfrac12\ \text{mod}\ 1
\end{array}\right.,\ \forall \mathrm{u}\in \T.
\]
The maps $F_{N,\epsilon}$ are all expanding piecewise affine maps with expanding rate $2(1-\epsilon)\in (1,2]$ (see Appendix \ref{A-EPWAMAPS} for a summary of related notions). Their atoms are determined by the pairwise distances between the coordinates $\mathrm{u}_i$, whether they are smaller or larger than $\tfrac12$. Moreover, the map $g$ commutes with the inversion symmetry $-\text{Id}|_{\T}$; likewise the maps $F_{N,\epsilon}$ commute with $-\text{Id}|_{\T^N}$. The $F_{N,\epsilon}$ also commute with every element of the {\bf group $\Pi_N$ of the permutations} of the coordinates $\{\mathrm{u}_i\}_{i=1}^N$.
Altogether, the $F_{N,\epsilon}$ can be considered as an elementary model of a system of $N$ chaotic units in interaction, where the discontinuities induced by $g$ play the role of nonlinearities. 

For every $N$, the map $F_{N,\epsilon}$ can be shown to have an ergodic acim when $\epsilon$ is small enough ({\sl viz.}\ expanding rate close to 2). Moreover, numerical simulations showed evidences of the breakdown of the inversion symmetry in the long-term dynamics, when $\epsilon$ is close enough to $\tfrac12$ ({\sl ie.}\ expanding rate close to 1) \cite{F14,F20}.
The essential characteristics of this phenomenology is given in Appendix \ref{A-PHENO}, which in particular describes the systematic symmetric and asymmetric features of the various acim. This appendix also introduces an original representation of the trajectories that facilitates the visualization of these features. 

The numerical evidences of loss of ergodicity have been partly confirmed by analytic and/or computer assisted proofs of the emergence of an {\bf asymmetric acim}, namely an acim whose support is disjoint from its image under the coordinate sign inversion. That symmetry then implies the existence of a pair of acim with disjoint supports, which suffices to ensure that ergodicity in $F_{N,\epsilon}$ cannot hold. 
The analytic proofs applied to $N\in [3,4]$, see \cite{F14,FS22,S18,SB16}, the computer-assisted ones to $N\in [3,6]$, see \cite{F20}. In both cases, they consisted in proving the existence of asymmetric invariant unions of polytopes ({\bf AsIUP}), namely invariant unions of polytopes ({\bf IUP}) that are disjoint from their image under the symmetry (see again Appendix \ref{A-EPWAMAPS} for the definitions). 

While the proofs have been designed to be deployed in arbitrary dimension, the details of the dynamics of the AsIUP are specific to the value of $N$ under consideration. No simple mechanism has emerged that could be naturally extended to an arbitrary value of $N$ (NB: Specific limitations to such an extension are discussed at the beginning of Section \ref{S-MULTACIM} below). 
Accordingly, this paper more modestly aims to provide instances of families of maps in arbitrary dimension that capture some characteristics of the $F_{N,\epsilon}$ while exhibiting provable emergence of AsIUP via a simple systematic mechanism when their expanding rate is close to 1. 

\subsection{Presentation of the results} 
\label{sec:presentationresults}
\begin{figure}[ht]
\begin{center}
\includegraphics*[width=130mm]{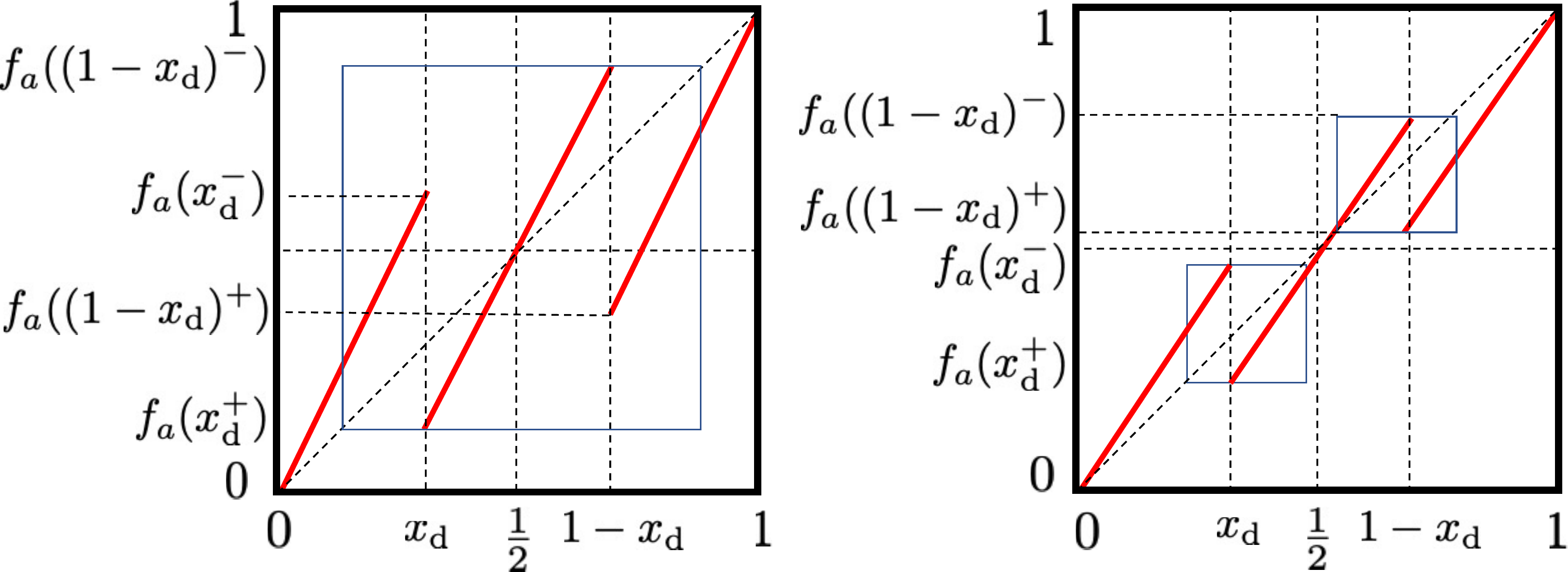}
\end{center}
\caption{Symmetry-breaking loss of ergodicity in the family $\{f_a\}_{a\in (1,2)}$ of symmetric piecewise affine Lorenz-type maps with three branches. {\sl Left.}\ When the slope $a$ is close to 2, the map $f_a$ is ergodic with unique acim. {\sl Right.}\ When $a$ is close to 1, the map $f_a$ has two acim with disjoint supports.}
\label{LORENZMAPS}
\end{figure}

A natural source of inspiration for our study is to consider a simple example in one dimension, namely the family $\{f_a\}_{a\in (1,2)}$ of Lorenz-type maps with three branches, see Figure \ref{LORENZMAPS} and Appendix \ref{A-LORENZ}. As the $F_{N,\epsilon}$, the maps $f_a$ are expanding piecewise affine and they commute with an inversion symmetry. Moreover, in agreement with the $F_{N,\epsilon}$ phenomenology described above (recall that the expanding rate of $F_{N,\epsilon}$ is equal to $2(1-\epsilon)$), $f_a$ is ergodic with unique acim when the expanding rate $a$ is close to 2, and has two acim with disjoint supports when $a$ is near 1.

Accordingly, we aim to provide some multi-dimensional analogues of the families of maps $\{f_a\}$ that capture some characteristics of the $F_{N,\epsilon}$. Focus will be made on the emergence of AsIUP. In particular, the multi-dimensional maps will be self-maps of the $(N-1)$-simplex $S_{N-1}$, where given an arbitrary $d\in\N$, $S_d$ is defined by 
\[
S_{d}=\left\{x=(x_i)_{i=1}^{d}\in \R_+^{d}\ :\ x_1+\cdots +x_{d}<1\right\}.
\]
Indeed, a preliminary analysis in Section \ref{S-TORUSMAPS} concludes that AsIUP of the coupled map $F_{N,\epsilon}$ in $\T^N$ can be deduced from AsIUP of some related projected map in $S_{N-1}$, denoted $G_{N-1,\epsilon}$.\footnote{Notice that similar $N-1$-dimensional reductions of the coupled map $F_{N,\epsilon}$ dynamics have already been identified \cite{F14,SB16}. Yet these reductions yielded maps of $\T^{N-1}$ or of $[0,1]^{N-1}$, not of $S_{N-1}$.} More precisely, an elementary reduction theory is developed in Section \ref{S-PROJECT}, for arbitrary maps of $\T^N$ that commute with every permutation of given $N-1$ coordinates. In the piecewise affine case, the theory shows that an IUP of the original map can be ensured by an IUP of some projected map of the subset $I_N\subset \R^N$ of points with increasing coordinates, see Lemma \ref{MULTSETS}. Considerations about the transfer of an additional symmetry, typically an inversion symmetry, to a symmetry of the projected map are given in Section \ref{S-ADDSYM}. For more specific mappings such as $F_{N,\epsilon}$, a further change of variables implies that it suffices to consider IUP of the map $G_{N-1,\epsilon}$, which also inherits some inversion symmetry from the $-\text{Id}|_{\T^N}$ of the $F_{N,\epsilon}$ (Section \ref{S-COUPLEDM}).

As a consequence, the desired multi-dimensional analogues of the map $f_a$ will be maps of $S_{N-1}$ that share some systematic features of the maps $G_{N-1,\epsilon}$. These features and the proposed mechanism for AsIUP are given in Section \ref{S-INSPIR}. In particular, Section \ref{S-OUTER} identifies the inversion symmetry, the outer atoms of the atomic collection and the dynamics in these atoms. This section also identifies the changes in the atomic collection as $N$ varies and hence, identifies some cause of the absence of full systematic in the features of the $G_{N-1,\epsilon}$. 

Inspiration for a systematic mechanism for AsIUP in the arbitrary $d$-simplex is obtained in Section \ref{S-G2} from a thorough analysis of the dynamics of the family $\{G_{2,\epsilon}\}_{\epsilon\in (0,\frac12)}$, which turns out to show the same ergodic/non-ergodic features as the family $\{f_a\}$. As for the $f_a$, the mechanism consists in ensuring, for sufficiently small expanding rate, the existence of a simply-connected IUP lying across $A\cup B$ where $A$ is an outermost atom and the atom $B$ is adjacent to $A$ (see Claim \ref{FEAT2D} and Fig.\ \ref{PARTITION2SIMPLEX}). As for the $f_a$, the map restriction to $A$ (resp.\ to $B$) expands away from the corresponding fixed point so that the points are eventually mapped into $B$ (resp.\ into $A$). In $B$, this trend is combined with the action of a permutation in some basis attached to the fixed point. The latter feature is a purely multi-dimensional characteristic that has no analogue in dimension 1. In addition, the dynamics in the atoms $A$ and $B$ are such that the IUP does not intersect its symmetric image; hence providing an AsIUP.
 
 With all the necessary elements being identified, the main result of the paper can finally be stated and proved (Theorem \ref{MAINRESULT} in Section \ref{S-MULTACIM}). This statement claims the existence, for an arbitrary integer $d\geq 3$, of a family of inversion-symmetric maps of the simplex $S_d$ that reproduce the common features of the maps $G_{d,\epsilon}$ and have AsIUP when the expanding rate is close enough to 1.  The proof essentially consists in constructing the maps so that they exhibit some multi-dimensional extension of the symmetry-breaking mechanism in the family $\{G_{2,\epsilon}\}_{\epsilon\in (0,\frac12)}$. In particular, focus is made on constructing the analogue of the restriction to the atom $B$ above. The construction and associated mechanism are based on simple principles which should extend to a broad range of maps with permutation and inversion symmetries. 

\section{Projection procedure for maps of the torus with permutation symmetries}\label{S-TORUSMAPS} 
This section introduces a projection procedure for maps of the torus $\T^N$ which commute with every element of the {\bf group $\Pi_{N-1}$ of the permutations} of the first $N-1$ coordinates of $\mathrm{u}\in\T^N$.\footnote{The projection procedure does not need that the map commutes with every element of $\Pi_N$, nor gain any benefit from that assumption, see the discussion in Section \ref{S-ADDSYM}.} Aiming at reducing every orbit generated by this symmetry group to a single point in phase space, the projection maps the torus to a subset $I_N$ of $\R^N$ of points with increasing coordinates, and subsequently to $S_{N-1}\times [0,1)$ by conjugacy. This procedure is particularly relevant in the case of invariant sets that consist of the orbit under $\Pi_{N-1}$ of a single connected component because these sets become simply connected invariant sets in the reduced dynamics. Instance of such sets have been observed in the phenomenology of the maps $F_{N,\epsilon}$ (see Appendix \ref{A-PHENO}). 

Furthermore, the procedure ensures that the existence of disjoint invariant sets/IUP of the induced map implies the same property for the original map in $\T^N$. In addition, conditions will be identified for an additional symmetry of the original map - typically an inversion symmetry - to be transferred to the projected map. In this setting, natural candidates for disjoint invariant sets of the projected map will be AsIUP. 

In the case of the maps $F_{N,\epsilon}$, their particular form implies that the corresponding map of $S_{N-1}\times [0,1)$ is a skew-product dynamical system whose base map $G_{N-1,\epsilon}$ acts in $S_{N-1}$. Moreover, the inversion symmetry $-\text{Id}|_{\T^N}$ is shown to transfer to the base map in $S_{N-1}$, see Fig.\ \ref{REDUCDIAGRAM} below for an illustration of the whole procedure in this case. Altogether, it suffices to prove the existence of AsIUP for the base map $G_{N-1,\epsilon}$ in order to conclude the existence of two acim with disjoint supports for $F_{N,\epsilon}$.

Actually, the same reduction to a map $G_{\rho,\epsilon}$ of $S_{N-1}$ applies to the coupled maps $F_{\rho,\epsilon}$ (see Appendix \ref{A-CLUSTERMAPS}) with distribution $\rho=(\rho_i)_{i=1}^N$ where $\rho_1=\rho_2=\cdots =\rho_{N-1}$. When $\rho_N$ differs from the other $\rho_i$, these maps $F_{\rho,\epsilon}$ only commute with the permutations in $\Pi_{N-1}$ (and not with the permutations in $\Pi_N\setminus \Pi_{N-1}$). In the main text, we keep considering $F_{N,\epsilon}$ for simplicity. We refer to Appendix \ref{A-CLUSTERMAPS} for those features that are specific to the more general $F_{\rho,\epsilon}$.

\subsection{The projection procedure and its consequences for the multiplicity of invariant sets}\label{S-PROJECT}
In order to define the projection procedure, we need to introduce various basic notions associated with the dynamics in $\T^N$ and its permutation symmetries. Given $N\in\N$, let
\[
\T_\ast^N=\left\{\mathrm{u}\in \T^N\ :\ \mathrm{u}_i\neq \mathrm{u}_j\ \text{mod}\ 1,\ \forall i\neq j\in [1,N]\right\}.
\]
be the set of elements $\mathrm{u}\in\T^N$ whose coordinates are all distinct. The reason for considering $\T_\ast^N$ instead of $\T^N$ will be given below. Let then the map $P$ be defined by
 \[
(P \mathrm{u})_i=\mathrm{u}_i+\lfloor \mathrm{u}_N-\mathrm{u}_i\rfloor-\lfloor \mathrm{u}_N\rfloor,\quad \forall i\in [1,N].
\]
This map is well-defined as a one-to-one mapping from $\T^N$ into $\R^N$. The set $D_\ast^N=P\T_\ast^N$ is a fundamental domain of $\T_\ast^N$ (namely, every element of $\T_\ast^N$ can be represented by a unique element in $D_\ast^N$), which reads 
\[
D_\ast^N=\left\{u\in \R^{N-1}\times [0,1)\ :\ u_i- u_j\in\R\setminus\Z,\ \forall i\neq j\in [1,N]\ \text{and}\ 0<u_N-u_i<1,\ \forall i\in [1,N-1]\right\}.
\]
Let also $I_N\subset D_\ast^N$ be the subset of points with increasing coordinates, namely 
\[
I_N =\left\{u\in \R^{N-1}\times [0,1)\ :\ u_1<u_2<\cdots <u_{N-1}<u_N<u_1+1\right\}.
\]
For the sake of notations, we use the same symbol $\pi$ for a transformation that permutes the first $N-1$ coordinates in $\T^N$ and in $\R^N$ respectively (and likewise for $\Pi_{N-1}$). Given $u\in D_\ast^N$, let $\pi_u\in \Pi_{N-1}$ be such that $\pi_u u\in I_N$. The reason for dealing with $\T_\ast^N$ instead of $\T^N$ is that the ordering permutation $\pi_u$ is unique when $u\in D_\ast^N$. Notice also that, for every $u\in D_\ast^N$, the map $v\mapsto \pi_uv$ is invertible on $\R^N$. Moreover, we obviously have $\pi I_N\subset D_\ast^N$ for every $\pi\in\Pi_{N-1}$.
\begin{Claim}
Given a map $F:\T^N\circlearrowleft$ which commutes with every $\pi\in \Pi_{N-1}$, let $\mathrm{F}:\R^N\circlearrowleft$ be defined by $\mathrm{F}=P\circ F\circ P^{-1}$. The restriction $\mathrm{F}|_{D_\ast^N}$ is entirely determined by its action on $I_N$, viz.\ we have
\begin{equation}
\mathrm{F} u=\pi_u^{-1}\circ \mathrm{F}|_{I_N}\circ \pi_u u,\ \forall u\in D_\ast^N.
\label{REDUC}
\end{equation}
\end{Claim}
\noindent
{\sl Proof:} The map $P$ and its inverse commute with every transformation in $\Pi_{N-1}$; hence so does $\mathrm{F}$ by the assumption on $F$. We then have 
\[
\mathrm{F}|_{I_N} v=\mathrm{F} v = \pi_u\circ \pi_u^{-1} \circ \mathrm{F} v =\pi_u\circ \mathrm{F}\circ \pi_u^{-1} v,\ \forall v\in I_N,u\in D_\ast^N.
\]
In particular, for $v=\pi_u u$, we get $\mathrm{F}|_{I_N}\circ \pi_u u=\pi_u\circ \mathrm{F}u$ from where the relation \eqref{REDUC} immediately follows. \hfill $\Box$

Assuming in addition that the map $F$ is non-singular, namely that the pre-images of zero Lebesgue measure sets have zero Lebesgue measure, so that $\mathrm{F}$ is also non-singular and then the complement set $I_N\setminus \left(I_N\cap \mathrm{F}^{-1}D_\ast^N\right)$ has zero Lebsegue measure. The relation \eqref{REDUC} suggests to consider the {\bf projected map} ${\cal F}$ defined in $I_N\cap \mathrm{F}^{-1}D_\ast^N$ by 
\[
u\mapsto {\cal F}u=\pi_{\mathrm{F} u}\circ \mathrm{F}|_{I_N} u.
\] 
The important features of ${\cal F}$ for our purpose, especially in the piecewise affine case, are identified in the following statement.
\begin{Lem}
(i) Assume that ${\cal F}$ has two disjoint forward invariant sets in $I_N\cap \mathrm{F}^{-1}D_\ast^N$. Then, the same property holds for $\mathrm{F}$ in $D_\ast^N$, and hence for $F$ in $\T_\ast^N$.

\noindent
(ii) If $F$ is a non-singular piecewise affine map, then there exist atomic collections in $D_\ast^N$ and in $I_N$ respectively so that the maps $\mathrm{F}$ and ${\cal F}$ are non-singular piecewise affine maps. Moreover, assume that ${\cal F}$ has two disjoint IUP in $I_N\cap \mathrm{F}^{-1}D_\ast^N$. Then, the same property holds for $\mathrm{F}$ in $D_\ast^N\cap \mathrm{F}^{-1}D_\ast^N$, and hence for $F$ in $\T_\ast^N\cap F^{-1}\T_\ast^N$.
\label{MULTSETS}
\end{Lem}
\noindent
Naturally, the converse statement cannot be true because of the equality 
\[
\pi_u=\pi_u\circ \pi,\ \forall \pi\in \Pi_{N-1},
\]
implies that every trajectory $\{F^t\mathrm{u}\}_{t\in\N}$ of $F$ and its image trajectory $\{\pi\circ F^t\mathrm{u}\}_{t\in\N}$ are mapped onto the same trajectory of ${\cal F}$. Yet, Lemma \ref{MULTSETS} can serve to detect distinct invariant sets/IUP of $F$ that consist of distinct orbits of the symmetry group $\Pi_{N-1}$.   
\smallskip

\noindent
{\sl Proof of the Lemma.} {\em (i)} Assume that $A,B\subset I_N\cap \mathrm{F}^{-1}D_\ast^N$ with $A\cap B\neq 0$ are two invariant sets of ${\cal F}$. Then both union sets $\bigcup_{\pi\in \Pi_{N-1}}\pi A,\bigcup_{\pi\in \Pi_{N-1}}\pi B\subset D_\ast^N$ must be disjoint invariant sets of $\mathrm{F}$.

To see this, assume that $u\in \bigcup_{\pi\in \Pi_{N-1}}\pi A$. Then $u\in \pi A\subset D_\ast^N$ for some $\pi\in\Pi_{N-1}$. Also the symmetry of $F$ implies that $\pi A\subset \mathrm{F}^{-1}D_\ast^N$ for every $\pi\in\Pi_{N-1}$, which in particular yields $\pi_uu\in I_N\cap \mathrm{F}^{-1}D_\ast^N$. Relation \eqref{REDUC} then implies
\[
\mathrm{F} u=\pi_u^{-1}\circ \mathrm{F}|_{I_N}\circ \pi_u u=\pi_u^{-1}\circ \pi_{\mathrm{F}\circ \pi_u u}^{-1}\circ \pi_{\mathrm{F}\circ \pi_u u}\circ \mathrm{F}|_{I_N}\circ \pi_u u=\pi_u^{-1}\circ \pi_{\mathrm{F}\circ \pi_u u}^{-1}v 
\]
where $v=\pi_{\mathrm{F}\circ \pi_u u}\circ \mathrm{F}|_{I_N}\circ \pi_u u={\cal F}\circ \pi_u u\in A$ because $\pi_u u\in A$ and $A$ is invariant under ${\cal F}$. In other terms, $\mathrm{F} u\in \pi' A$ for some $\pi'\in \Pi_{N-1}$, proving invariance. 

Moreover, that the union sets $\bigcup_{\pi\in \Pi_{N-1}}\pi A$ and $\bigcup_{\pi\in \Pi_{N-1}}\pi B$ are disjoint is also immediate. Firstly, when $\pi\neq \pi'$, we must have $\pi A\cap \pi' B=\emptyset$ because these sets belong to distinct regions of $D_\ast^N$ (distinct relative ordering of the coordinates). Secondly, if we had $\pi A\cap \pi B\neq\emptyset$ for some $\pi$, then we would have $A\cap B\neq \emptyset$ since $\pi$ is one-to-one, which contradicts the initial assumption. 

\noindent
{\em (ii)} Assume that $F$ is a non-singular piecewise affine map and let $\{A_\omega\}$ be its atomic collection (see Appendix \ref{A-EPWAMAPS}). Then $\mathrm{F}$ is also a non-singular piecewise affine map for the atomic collection, say $\{A'_{\omega'}\}$, defined by the refinement of the image collection $\{P A_\omega\}$ by the level sets of the functions $\{\lfloor (Fu)_N-(Fu)_i\rfloor-\lfloor (Fu)_N\rfloor\}_{i\in [1,N]}$. The conjugacy $P$ implies that any IUP $\bigcup_k P_k$ of $\mathrm{F}$ induces an IUP $\bigcup_k P^{-1}P_k$ of $F$ for the refined atomic collection $\{P^{-1}A'_{\omega'}\}$. Moreover, two distinct IUP of $\mathrm{F}$ induce two distinct IUP of $F$. 

In addition, the permutation symmetry group implies that for every $\omega'$, we have $A'_{\omega'}=\pi A'_{\omega'_+}$ where $\pi\in \Pi_{N-1}$ and the index $\omega'_+$ are such that $A'_{\omega'_+}\subset I_N$. Then, the map ${\cal F}$ is a piecewise affine map for the atomic collection $\{A''_{\omega''}\}$ defined by the refinement of $\{A'_{\omega'_+}\}$ by the sets in $I_N$ in which the ordering of the coordinates $((\mathrm{F} u)_i)_{i=1}^{N-1}$ is constant (namely those sets in which the permutation $\pi_{\mathrm{F} u}$ does not depend on $u$). 

Now, a similar reasoning as in the proof of {\em (i)} shows that if $\bigcup_kP_k$ and $\bigcup_{k'}P'_{k'}$ are disjoint IUP for ${\cal F}$, then $\bigcup_kP_k$ and $\bigcup_{k'}P'_{k'}$ must be disjoint IUP for $\mathrm{F}$. \hfill $\Box$
   
\subsection{Maps with additional symmetries}\label{S-ADDSYM}
Following Lemma \ref{MULTSETS}, a natural setting for the existence of multiple invariant sets/IUP for the projected map ${\cal F}$ is when this map has some symmetry, so that one can investigate the existence of asymmetric invariant sets/AsIUP. Accordingly, we need to determine those conditions that ensure that a (additional) symmetry of the original map $F$ transfers to one for ${\cal F}$. This is precisely the purpose of the following statement.
\begin{Lem}
Assume that a transformation $S:\T_\ast^N\circlearrowleft$ commutes with $F$ and that the induced transformation $\Sigma=P\circ S\circ P^{-1}$ on $D_\ast^N$ has a proper representation on $I_N$, ie.\ there exists $\sigma_\Sigma:I_N\circlearrowleft$ such that 
\[
\sigma_\Sigma \circ \pi_u u=\pi_{\Sigma u}\circ \Sigma u,\ \forall u\in D_\ast^N.
\]
Then, ${\cal F}$ commutes with $\sigma_\Sigma$ on $I_N\cap F^{-1}\T_\ast^N$. 
\label{REDUCSYM}
\end{Lem}
\noindent
{\sl Proof.} Throughout the proof we use the symbol $\sigma$ to denote $\sigma_\Sigma$. For every $u\in D_\ast^N\cap \mathrm{F}^{-1}D_\ast^N$, we have $ \pi_{\mathrm{F} u }\circ \mathrm{F} u\in I_N\subset D_\ast^N$; hence using the characterization of $\sigma$ above and $\Sigma\circ \mathrm{F}=\mathrm{F}\circ \Sigma$, we get
\[
\sigma\circ \pi_{\mathrm{F} u }\circ \mathrm{F} u= \pi_{\Sigma \circ \mathrm{F} u}\circ \Sigma\circ \mathrm{F} u=\pi_{\mathrm{F}\circ \Sigma u}\circ \mathrm{F}\circ \Sigma u.
\]
On the other hand, we have $\pi_u=\text{\rm Id}$ on $I_N$ and then for $u\in I_N$
\begin{align*}
\pi_{\mathrm{F}\circ \sigma u}\circ \mathrm{F}\circ \sigma u&=\pi_{\mathrm{F}\circ \pi_{\Sigma u}\circ \Sigma u}\circ \mathrm{F}\circ \pi_{\Sigma u}\circ \Sigma u=\pi_{\pi_{\Sigma u}\circ \mathrm{F}\circ  \Sigma u}\circ  \pi_{\Sigma u}\circ \mathrm{F}\circ \Sigma u\\
&=\pi_{\mathrm{F}\circ \Sigma u}\circ \mathrm{F}\circ \Sigma u
\end{align*}
where the second equality follows from the fact that $\mathrm{F}$ commutes with $\pi_{\Sigma u}$ and the second line follows from the fact that 
\begin{equation}
\pi_{\pi' u}\circ \pi' u=\pi_u u,\ \forall \pi'\in \Pi_{N-1}.
\label{DEGEPI}
\end{equation}
\hfill $\Box$

Since $F$ commutes with every $\pi\in \Pi_{N-1}$, it follows that,  when $S$ satisfies the conditions of Lemma \ref{REDUCSYM}, every transformation $\pi\circ S$ induces a transformation $\pi\circ \Sigma$ with proper representation on $I_N$. Yet, that representation is identical to that of the original symmetry $S$, as our next claim states.
\begin{Claim}
We have $\sigma_\Sigma=\sigma_{\pi\circ \Sigma}$ for every $\pi\in \Pi_{N-1}$.
\end{Claim}
\noindent
 {\sl Proof.} This is immediate from the relation \eqref{DEGEPI} and the fact that $\pi_{\pi' u} u=\pi_u u$ for every $\pi'\in \Pi_{N-1}$.\hfill $\Box$
 
 \begin{Exa}
 The inversion of coordinate signs $S=-\text{\rm Id}|_{\T^N}$ by 
\[
(S \mathrm{u})_i=-\mathrm{u}_i\ \text{mod}\ 1,\ \forall i\in [1,N],
\]
induces the transformation $\Sigma=P\circ S\circ P^{-1}$ on $D_\ast^N$ whose explicit expression reads (after simple algebra)
\[
(\Sigma u)_i=\delta_{i,N}-\delta_{u_N,0}-u_i,\ \forall i\in [1,N].
\]
Clearly, $\Sigma$ has a proper representation $\sigma_\Sigma$ on $I_N$ given by   
\begin{equation}
(\sigma_\Sigma u)_i=\left\{\begin{array}{ccl}
-\delta_{u_N,0}-u_{N-i}&\text{if}&i\in [1,N-1]\\
1-\delta_{u_N,0}-u_N&\text{if}&i=N
 \end{array}\right.
 \label{INVSYM}
\end{equation}
\end{Exa}
\begin{Rem}
The left cyclic permutation $K$ of the $N$ coordinates in $\T^N$ defined by
\[
 (K \mathrm{u})_i=\left\{\begin{array}{ccl}
 \mathrm{u}_{i+1}\ \text{mod}\ 1&\text{if}&i\in [1,N-1]\\
 \mathrm{u}_1\ \text{mod}\ 1&\text{if}&i=N
 \end{array}\right.
 \]
 induces the following transformation $\kappa=P\circ K\circ P^{-1}$ on $D_\ast^N$
 \[
 (\kappa u)_i=\left\{\begin{array}{ccl}
u_{i+1}+\lfloor u_1-u_{i+1}\rfloor-\lfloor u_1\rfloor&\text{if}&i\in [1,N-1]\\
 u_1-\lfloor u_1\rfloor&\text{if}&i=N
 \end{array}\right.
 \]
 This map has no proper representation in $I_N$. Indeed, for every $\pi\in \Pi_{N-1}$ that affects the first coordinate, we have 
 \[
 (\pi_{\kappa \circ \pi u}\circ \kappa \circ \pi u)_N=(\kappa \circ \pi u)_N\neq (\kappa u)_N=(\pi_{\kappa u}\circ \kappa \circ \pi u)_N,\ \forall u\in D_\ast^N.
 \]
 Yet, we have $\pi_{\pi u}\circ \pi u=\pi_u u$; hence the equality in Lemma \ref{REDUCSYM} cannot hold. 
\end{Rem}
\noindent
This remark shows that the cyclic permutation symmetry cannot transfer to the projected map ${\cal F}$. In order words, that in addition to $\Pi_{N-1}$, the original map $F$ also commutes with $K$ (and hence with every permutation of the $N$ coordinates, by composition) does not bring any additional symmetry to ${\cal F}$. 

However when $F$ commutes with every permutation in $\Pi_N$, or more generally, when it commutes (only) with every element of the group $\Pi_{i_1,\cdots ,i_{N-1}}$ of the permutations of the $(N-1)$ coordinates indexed by $\{i_1,\cdots ,i_{N-1}\}$, a similar projection procedure as in the previous section can be defined, which is adapted to the $(N-1)$-uple under consideration. Naturally, the fundamental domain and corresponding projection $P$ depend on this $(N-1)$-uple, as well as do the representation of $F$ on the corresponding set $I_N$ of points with increasing coordinates and that of the inversion of sign coordinates $S$.

In short terms, when $F$ commutes with every permutation of all coordinates, both the projected map and the representation of the inversion of coordinate signs on $I_N$ are not unique and depend on the choice of the fundamental domain, see Appendix \ref{A-ALTERNAT} for examples.

\subsection{Application to the coupled maps}\label{S-COUPLEDM}
For the sake of notation, let $d=N-1$. As already pointed out, the coupled maps $F_{N,\epsilon}$ commute with every $\pi\in\Pi_{d}$, and also with the inversion symmetry $S=-\text{Id}|_{\T^N}$. The arguments above imply that when the corresponding projected map ${\cal F}_{N,\epsilon}$ defined by 
\[
u\mapsto {\cal F}_{N,\epsilon}u=\pi_{\mathrm{F}_{N,\epsilon} u}\circ \mathrm{F}_{N,\epsilon}|_{I_N} u,\ u\in I_N\cap \mathrm{F}_{N,\epsilon}^{-1}D_\ast^N,
\] 
(where $\mathrm{F}_{N,\epsilon}=P\circ F_{N,\epsilon}\circ P^{-1}$) has an AsIUP with respect to the inversion symmetry $\sigma_\Sigma$ defined by \eqref{INVSYM}, then $F_{N,\epsilon}$ must have two acim with disjoint supports.

In addition, the specific form of the expression of $F_{N,\epsilon}$, namely that it consists of the sum of a multiple of the identity on $\T^N$ and a map that only depends on the coordinates differences $\mathrm{u}_j-\mathrm{u}_i$, implies a further reduction. To see this, let $\phi_N$ be defined by 
\begin{equation*}
 (\phi_N u)_i=\left\{\begin{array}{ccl}
 u_{i+1}-u_i&\text{if}&i\in [1,d]\\
 u_N&\text{if}&i=N
 \end{array}\right.,\quad u\in\R^N.
 \end{equation*}
This map is one-to-one and we have $\phi_N I_N=S_{d}\times [0,1)$, where $S_d$ is the $d$-simplex introduced in Section \ref{sec:presentationresults}. Moreover, explicit computations yield the following statement.
\begin{Claim}
The conjugated map $\phi_N\circ {\cal F}_{N\epsilon}\circ \phi_N^{-1}$ is a skew-product dynamical system on $S_{d}\times [0,1)$, whose base map, say $G_{d,\epsilon}$, is a piecewise affine map from $S_{d}$ into itself. 
\label{CONJUG}
\end{Claim}
Moreover, we have $\phi_N\circ \sigma_\Sigma\circ \phi_N^{-1}=\sigma_{d}\times\sigma'_{1}$, where the inversion symmetries $\sigma_{d}$ and $\sigma'_{1}$ respectively act on $S_{d}$ and $[0,1)$, and are given by
\begin{equation}
(\sigma_{d} x)_i=\left\{\begin{array}{ccl}
x_{d-i}&\text{if}&i\in [1,d-1]\\
1-(x_1+\cdots+x_{d})&\text{if}&i=d
\end{array}\right.,\ x\in S_{d},
 \label{CHANGVAR}
\end{equation}
and
\[
\sigma'_{1} x= 1-\delta_{x,0}-x,\ x\in [0,1)
\]
As a consequence, all the maps $G_{d,\epsilon}$ commute with $\sigma_{d}$.

Now, one can show that if $\bigcup_kP_k$ is an AsIUP of $G_{d,\epsilon}$ with respect to $\sigma_{d}$, then $\phi_N^{-1}\left(\bigcup_k P_k\times [0,1)\right)$ is an AsIUP of the projected map ${\cal F}_{N,\epsilon}$ with respect to $\sigma_\Sigma$. Accordingly, it suffices to obtain an AsIUP in $S_{d}$ of $G_{d,\epsilon}$ in order to show the existence of two acim with disjoint supports in $\T^N$ for the original coupled maps $F_{N,\epsilon}$.

\begin{figure}[ht]
\begin{center}
\includegraphics*[width=160mm]{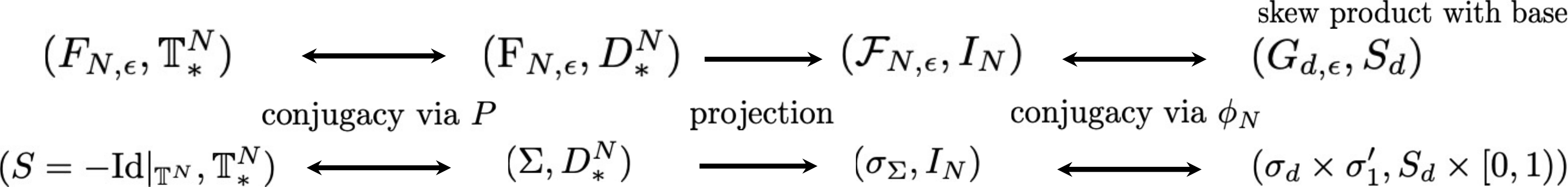}
\end{center}
\caption{Schematic representation of the whole reduction procedure for the coupled maps $F_{N,\epsilon}$ and their inversion symmetry $S=-\text{Id}|_{\T^N}$. The original system $(F_{N,\epsilon},\T_\ast^N)$ is first conjugated to $(\mathrm{F}_{N,\epsilon},D_\ast^N)$, then projected to $({\cal F}_{N,\epsilon},I_N)$, which is in turn conjugated to a skew-product system whose base is $(G_{d,\epsilon},S_d)$. Similar operations are applied to the symmetry $S$, which yield the symmetry $\sigma_d$ defined in \eqref{CHANGVAR} for the system $(G_{d,\epsilon},S_d)$.}
\label{REDUCDIAGRAM}
\end{figure}

A schematic summary of the whole reduction procedure associated with $F_{N,\epsilon}$ is given in Fig.\ \ref{REDUCDIAGRAM}.

\section{Inspiring features of the maps $G_{d,\epsilon}$}\label{S-INSPIR}
In this section, we identify some basic features of the reduced maps $G_{d,\epsilon}:S_d\circlearrowleft$ that will inspire the construction to come.
Focus is made on the outermost atoms of the atomic collection, namely those atoms that consist of $d$-simplexes whose facets are included in the facets of $S_d$ itself, and their adjacent atoms that are separated by a co-dimension 1 facet contained in the interior of $S_d$.\footnote{As we shall see below, the latter are genuine atoms only when $d\in [2,3]$. Otherwise, the map $G_{d,\epsilon}$ is not continuous on these sets. Moreover, its discontinuities depend both on $d$ and $\epsilon$. This is a cause of the absence of a simple systematic mechanism for the loss of ergodicity in the $F_{N,\epsilon}$ when $N\geq 5$.} In addition, the mechanism responsible for the emergence of AsIUP in $G_{2,\epsilon}$ is thoroughly analyzed.

\subsection{Characterisation of $G_{d,\epsilon}$ in the outer atoms of $S_d$}\label{S-OUTER}
\begin{figure}[ht]
\begin{center}
\includegraphics*[width=130mm]{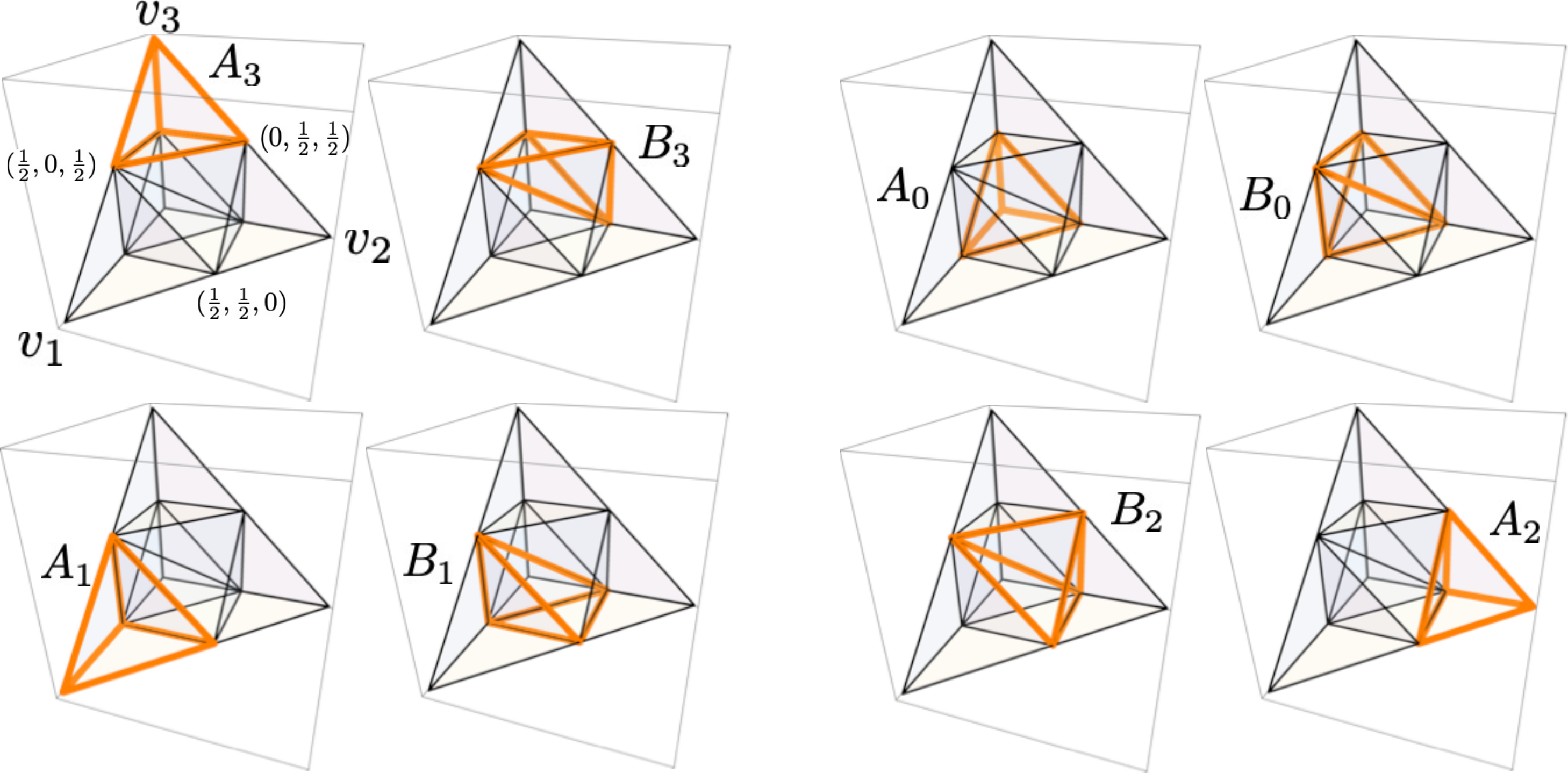}
\end{center}
\caption{The simplex $S_3$ and the composing simplexes $A_k$ and $B_k$. The collection $\{A_k,B_k\}_{k=0}^2$ is an atomic collection for $G_{3,\epsilon}$.}
\label{PARTITION3SIMPLEX}
\end{figure}

In order to state the features of $G_{d,\epsilon}$, we need to introduce and to describe the following subsets of $S_d$, see Fig.\ \ref{PARTITION3SIMPLEX} and the left panel in Fig.\ \ref{PARTITION2SIMPLEX}. 
\begin{itemize}
\item Given $d\in \N$ and $k\in [0,d]$, let 
\[
A_0=\left\{x\in S_{d}\ :\ x_1+\cdots +x_{d}<\tfrac12\right\}\quad\text{and}\quad A_k=\left\{x\in S_d\ :\ \tfrac12 <x_k\right\},\ k\in [1,d].
\]
\item For $d\geq 2$, let 
\begin{align*}
&B_0=\left\{x\in S_d\ :\ x_1+\cdots +x_{d-1},x_2+\cdots +x_{d}<\tfrac12 <x_1+\cdots +x_{d}\right\}\\
&B_1=\left\{x\in S_d\ :\ x_1,x_{2}+\cdots +x_{d}<\tfrac12<x_1+x_{2}\right\},\\
&B_k=\left\{x\in S_d\ :\ x_k<\tfrac12<x_{k-1}+x_k,x_k+x_{k+1}\right\},\ k\in [2,d-1],
\end{align*}
and 
\[
B_{d}= \left\{x\in S_d\ :\ x_1+\cdots +x_{d-1},x_{d}<\tfrac12<x_{d-1}+x_{d}\right\}.
\]
\end{itemize}
Let $\{v_k\}_{k=0}^d$ be the collection of the vertices of $S_d$, where $v_0$ is the origin and where the coordinates of the other $v_k$ satisfy $(v_k)_i=\delta_{i,k}$ for $i\in [1,d]$. 
\begin{Claim}
(i) The sets $A_k$ are pairwise disjoint $d$-simplexes included in $S_d$, whose vertices are $v_k$ and the middle points of the edges of $S_d$ issued from $v_k$.

\noindent
(ii) Each set $B_k$ is also a $d$-simplex included in $S_d$ and adjacent to $A_k$. The sets $B_k$ are pairwise disjoint if $d\geq 3$ and they all coincide for $d=2$. Moreover, we have $\bigcup_k A_k\cap \bigcup_kB_k=\emptyset$ and each $A_k\cup B_k$ forms a convex bipyramid.

\noindent
(iii) Recall the inversion symmetry $\sigma_d$ defined in \eqref{CHANGVAR}. We have $\sigma_dA_k=A_{d-k}$ for $k\in [0,\lceil\tfrac{d}2\rceil -1]$,\footnote{The symbol $\lceil\cdot\rceil$ denotes the ceiling function.} and if $d$ is even, we also have $\sigma_dA_{\frac{d}2}=A_{\frac{d}2}$. The same properties hold for the $B_k$.
\end{Claim}
In addition, notice that the Lebesgue measure of the complement set $S_d\setminus\bigcup_{k}(A_k\cup B_k)$ is zero for $d\in [2,3]$ and positive for $d\geq 4$. 

\noindent
{\sl Proof of the Claim.} {\em (i)} That $A_0$ (resp.\ $A_k$) is a $d$-simplex is a direct consequence of the fact that it can be obtained as the truncation of $S_d$ by the hyperplane $x_1+\cdots +x_d=\tfrac12$ (resp.\ $x_k=\tfrac12$). That the $A_k$  are pairwise disjoint is immediate from their definition and the constraints in $S_d$. 

\noindent
{\em (ii)} That $B_k$ is a simplex follows from the fact that is has $d+1$ facets that are given by the following independent inequalities
\begin{align*}
0<x_\ell,\ \ell\in [2,d-1]\ \text{and}\ x_1+\cdots +x_{d-1},x_2+\cdots +x_{d}<\tfrac12 <x_1+\cdots +x_{d}&\quad\text{if}\ k=0\\
0<x_\ell,\ \ell\in [3,d]\ \text{and}\ x_1, x_2+\cdots +x_{d}<\tfrac12<x_1+x_2&\quad\text{if}\ k=1&\\
0<x_\ell,\ \ell\in [1,d]\setminus [k-1,k+1],\ x_k<\tfrac12<x_{k-1}+x_k,x_k+x_{k+1}\ \text{and}\ x_1+\cdots +x_{d}<1&\quad\text{if}\ k\in [2,d-1]\\
0<x_\ell,\ \ell\in [1,d-2]\ \text{and}\ x_1+\cdots +x_{d-1}, x_{d}<\tfrac12<x_{d-1}+x_d&\quad\text{if}\ k=d
\end{align*}
That the $B_k$ are pairwise disjoint is immediate from their definition and the constraints in $S_d$. Moreover, the (only) facet of $A_k$ included in the interior of $S_d$ is also a facet of $B_k$; hence $A_k$ and $B_k$ must be adjacent sets. That $\overline{A_k\cup B_k}$ is convex is immediate from their definition. 

\noindent
{\em (iii)} Proved by direct computations. \hfill $\Box$

Now, the next statement describes the main properties of the restrictions $G_{d,\epsilon}|_{A_k}$ and $G_{d,\epsilon}|_{B_k}$.
\begin{Lem}
In addition to commuting with $\sigma_d$, the piecewise affine map $G_{d,\epsilon}$ has the following features for every $\epsilon\in \left(0,\tfrac12\right)$.

\noindent
(i) Every simplex $A_k$ is an atom of $G_{d,\epsilon}$ and the restrictions of $G_{d,\epsilon}$ to $A_0$ and to $A_k$ respectively write
\[
(G_{d,\epsilon}|_{A_0}x)_i=2(1-\epsilon) x_i\quad\text{\em and}\quad (G_{d,\epsilon}|_{A_k}x)_i=2(1-\epsilon) x_i+(2\epsilon-1)\delta_{i,k},\ i\in [1,d].
\]

\noindent
(ii) For any $k\in [0,d]$, the simplex $B_k$ is an atom of $G_{d,\epsilon}$ iff $d\in [2,3]$. 
\label{FEATURESCOUPLEDMAPS}
\end{Lem}
\noindent
This statement is a special case for uniform distributions $\rho$ ({\sl viz.}\ for $\varrho=\tfrac1{N}=\tfrac1{d+1}$) of Lemma \ref{A-FEATURESCOUPLEDMAPS} in Appendix \ref{A-CLUSTERMAPS}. We refer to that Appendix for a proof. 

In order words, Lemma \ref{FEATURESCOUPLEDMAPS} states that $\{A_k,B_k\}_{k=0}^d$ is an atomic collection of $G_{d,\epsilon}$ for $d\in [2,3]$, and also that $G_{d,\epsilon}$ has discontinuities inside every $B_k$ when $d\geq 4$. In the next section, we provide an analysis of the dynamics of $G_{2,\epsilon}$ and we establish the existence of AsIUP for $\epsilon$ near $\tfrac12$. 

\subsection{Analysis of the reduced map $G_{2,\epsilon}$}\label{S-G2}
As a preliminary comment to this section, we observe that for $d=1$, using that the interval $S_1=(0,1)=A_0\cup A_1\cup \{\tfrac12\}$, the first claim in Lemma \ref{FEATURESCOUPLEDMAPS} entirely determines (up to a set of zero Lebesgue measure) the one-dimensional map $G_{1,\epsilon}:S_1\circlearrowleft$.
This map is a particular case of a Lorenz map with two branches \cite{P79,W79} and has a unique ergodic acim for every value of $\varrho$ and $\epsilon$. 

\begin{figure}[ht]
\begin{center}
\includegraphics*[width=53mm]{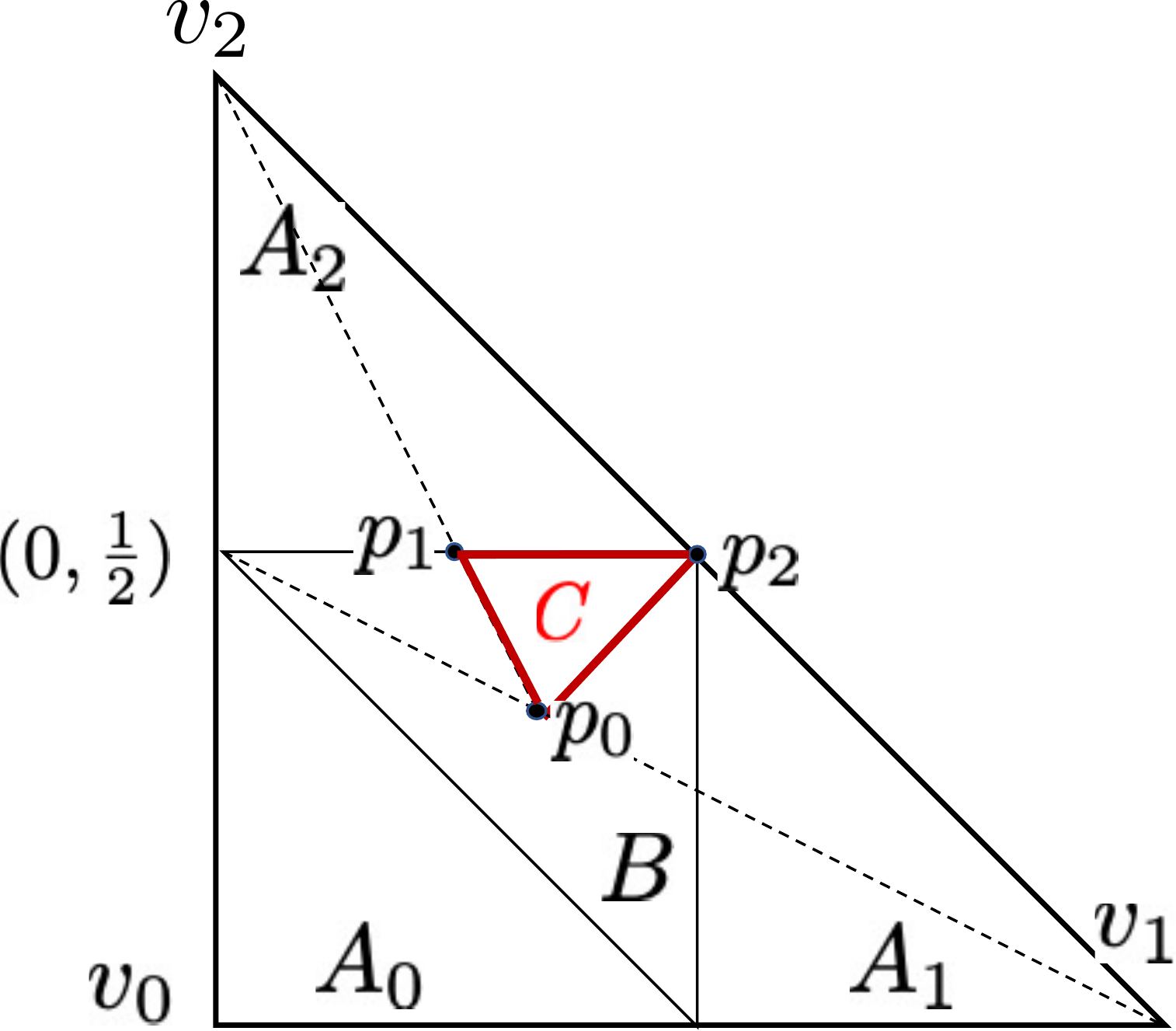}
\hspace{0.7cm}
\includegraphics*[width=47mm]{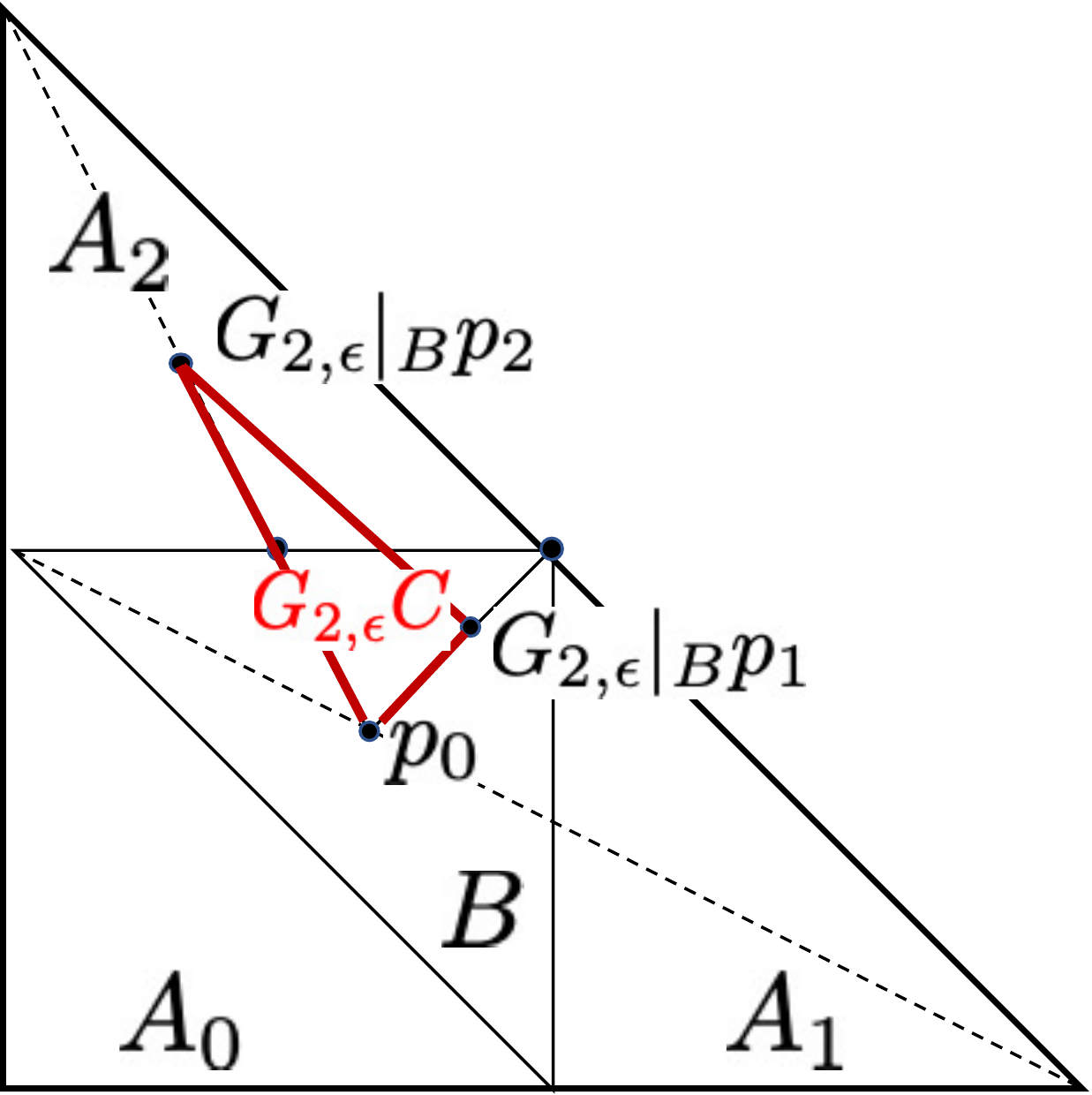}
\hspace{1cm}
\includegraphics*[width=42mm]{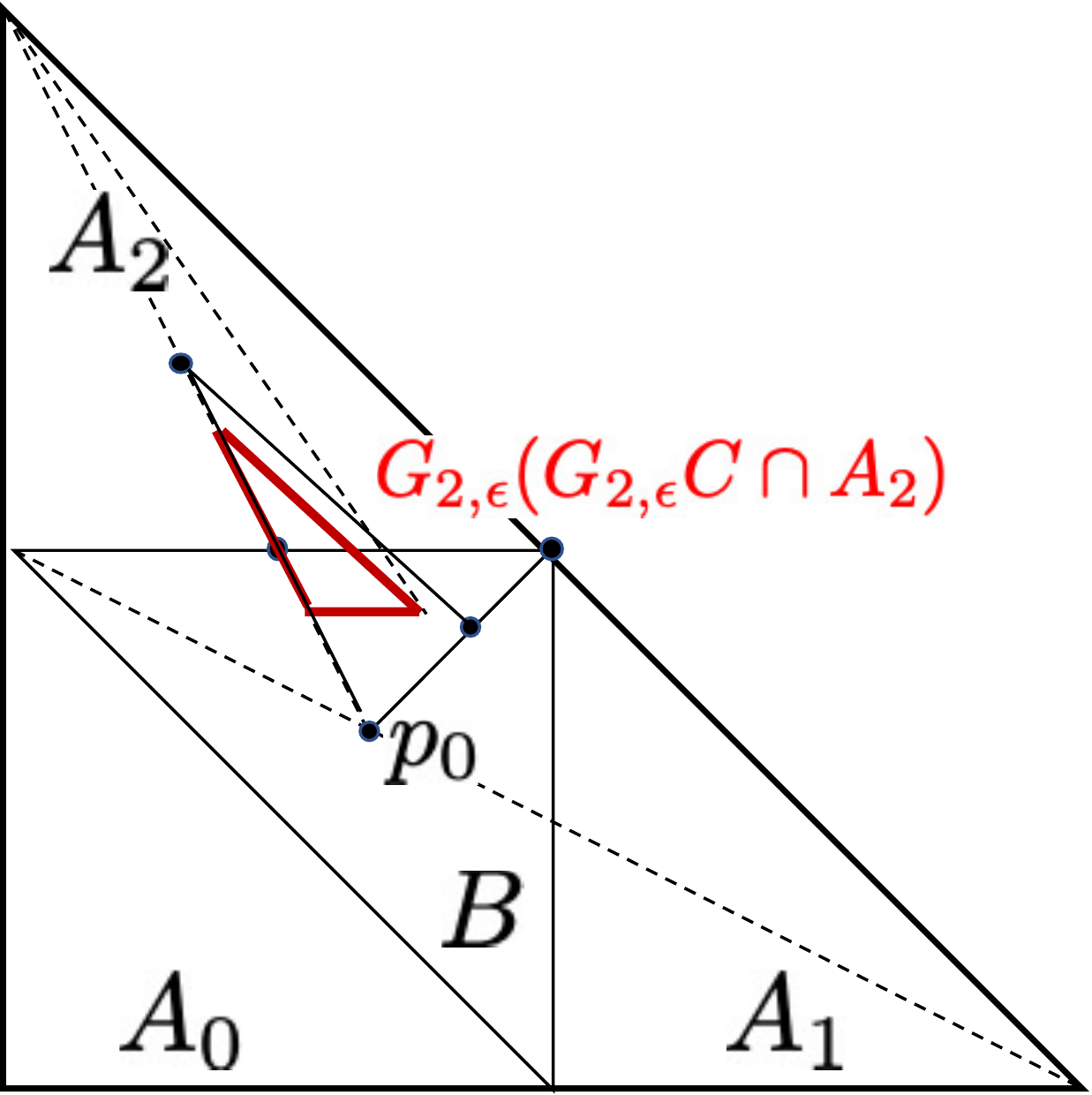}
\end{center}
\caption{{\sl Left.}\ The triangle $S_2$ and the composing triangles $A_0,A_1,A_2$ and $B:=B_0=B_1=B_2$ which form an atomic collection of $G_{2,\epsilon}$. The triangle $C$ has vertices $p_0,p_1,p_2$, which are defined in Claim \ref{FEAT2D}. {\sl Center.}\ The triangle $G_{2,\epsilon}C$ is obtained using the features listed in Claim \ref{FEAT2D}. {\sl Right.}\ The restriction $G_{2,\epsilon}|_{A_2}$ is an expanding affine map with  fixed point equal to $v_2$ and affine part given by $2(1-\epsilon)\mathrm{Id}$. Hence, when $\epsilon$ is close enough to 1, we have $G_{2,\epsilon}(G_{2,\epsilon}C\cap A_2)\subset G_{2,\epsilon}C$. In this case, the set $C\cup G_{2,\epsilon}C$ is an IUP of $G_{2,\epsilon}$.}
\label{PARTITION2SIMPLEX}
\end{figure}

For $d=2$, we have $S_2=\overline{A_0\cup A_1\cup A_2\cup B}$ where  $B:=B_0=B_1=B_2$ (see the left panel in Fig.\ \ref{PARTITION2SIMPLEX}). Lemma \ref{FEATURESCOUPLEDMAPS} states that the 2-dimensional map $G_{2,\epsilon}:S_2\circlearrowleft$ is an piecewise affine map with atoms $A_0,A_1,A_2$ and $B$. That statement also describes the action of the restrictions $G_{2,\epsilon}|_{A_k}$. As for the restriction $G_{2,\epsilon}|_{B}$, its expression is as follows (see equation \eqref{EXPRG} in Appendix \ref{A-FEATURESCOUPLEDMAPS})
\[
G_{2,\epsilon}|_Bx=(-2(1-\epsilon)x_1+1-\tfrac{2\epsilon}3,2(1-\epsilon)(x_1+x_2)+\tfrac{4\epsilon}3-1).
\]
An analysis of this expression (details not shown) reveals that $G_{2,\epsilon}|_{B}$ has the following characteristics.
\begin{Claim}
(i) The fixed point $p_0=\left(\tfrac13,\tfrac13\right)$ of $G_{2,\epsilon}|_{B}$ belongs to $B$.

\noindent
(ii) Let $p_1$ be the intersection point of the segment $[p_0v_2]$ and the edge $\overline{A_2}\cap\overline{B}$. Let $p_2$ be the intersection point (which exists) of the image segment $G_{2,\epsilon}|_B[p_0p_1]$ and $\overline{A_2}\cap\overline{B}$. In the basis formed by the vectors $p_0p_1$ and $p_0p_2$, the linear part of $G_{2,\epsilon}|_{B}$ is given by the following matrix
\[
2(1-\epsilon)\left(\begin{array}{cc}
0&\tfrac12\\
2&0\end{array}\right)
\]
\label{FEAT2D}
\end{Claim}
Let $C$ be the triangle with vertices $p_0,p_1$ and $p_2$. The properties of $G_{2,\epsilon}$ on $A_2\cup B$ imply that the set $C\cup G_{2,\epsilon}C\subset \overline{A_2\cup B}$ is an IUP of $G_{2,\epsilon}$ when $\epsilon$ is close enough to $\tfrac12$.\footnote{Notice that $G_{2,\epsilon}|_{B}$ is not expanding but its second iterate $(G_{2,\epsilon}|_{B})^2$ is. This is sufficient to conclude that the IUP must contains an acim.} The proof of this conclusion, which essentially consists in showing that $G_{2,\epsilon}(G_{2,\epsilon}C\cap A_2)\subset G_{2,\epsilon}C$ when the expansion rate $2(1-\epsilon)$ is sufficiently close to 1, is sketched on Figure \ref{PARTITION2SIMPLEX}. The proof itself is an adaptation {\sl mutatis mutandis} for $d=2$ of the proof of Proposition \ref{PROHd} below (details not shown).

Furthermore, recall that the map $G_{2,\epsilon}$ commutes with the transformation $\sigma_2$ defined by 
\[
\sigma_2(x_1,x_2)=(x_1,1-x_1-x_2)
\]
which is the reflection symmetry transverse to the line $(p_0v_1)$. Clearly, Figure \ref{PARTITION2SIMPLEX} shows that $C\cup G_{2,\epsilon}C$ is disjoint from its image under $\sigma_2$ (which lies in $\overline{A_0\cup B}$), and hence it is actually an AsIUP when it is an IUP.\footnote{Notice that the original map $F_{3,\epsilon}$ also commutes with every permutation of the $2$ coordinates $u_1$ and $u_3$. The projection procedure in that setting yields another map $G'_{2,\epsilon}$ of $S_2$ whose inversion symmetry is the reflection with respect to the diagonal $(x_1,x_2)\mapsto (x_2,x_1)$  (see Appendix \ref{A-ALTERNATSUB2}). That symmetry may appear more natural than the reflection with respect to $(p_0v_1)$ and the map $G'_{2,\epsilon}$ also has an AsIUP when $\epsilon$ is close to $\tfrac12$.} In particular, this simply connected AsIUP confirms the numerical phenomenology of the coupled map $F_{3,\epsilon}$ for $\epsilon$ close to $\frac12$ (see Appendix \ref{A-PHENO}).
 
As a side comment, notice that it can be proved that the map $G_{2,\epsilon}$ is locally eventually onto when $\epsilon$ is small enough (details not shown - see \cite{F14} for a proof for a similar symmetric map of the unit square), which implies uniqueness of the acim and hence ergodicity. In other words, the family of two dimensional maps $\{G_{2,\epsilon}\}_{\epsilon\in (0,\frac12)}$ has the same dynamical features as the family of one-dimensional Lorenz-type maps $\{f_a\}_{a\in (1,2)}$ given in the introduction.
 
\section{Symmetric maps of $S_d$ with multiple acim of disjoint supports}\label{S-MULTACIM}
Ideally, we would like to prove for $G_{d,\epsilon}$ with $d\geq 3$ arbitrary, the same emergence of multiple acim with disjoint (asymmetric) supports for $\epsilon$ near $\tfrac12$. 
More precisely, given the numerical phenomenology reported in Appendix \ref{A-PHENO}, we would like to prove the existence of simply connected AsIUP. 

However, the features of $G_{d,\epsilon}$ indicate that, if they exist, such AsIUP for $d\geq 3$ cannot be as simple as the one in $G_{2,\epsilon}$ described in the previous section. In particular, for $d=3$ the expression \eqref{EXPRG} in Appendix \ref{A-FEATURESCOUPLEDMAPS} of the image coordinate $(G_{3,\epsilon}|_{B_1}x)_3$ implies that no AsIUP can be included in $\overline{A_1\cup B_1}$ (nor in $\overline{A_3\cup B_3}$ by symmetry).\footnote{Actually, when $\epsilon$ is close to $\tfrac12$, $G_{3,\epsilon}$ has a sinply-connected AsIUP across $B_0$, $A_1$ and $B_1$ (details not shown).} In addition, one can check that the fixed point of $G_{3,\epsilon}|_{B_0}$ does not belong $S_3$; hence a feature as in statement {\sl (i)} of Claim \ref{FEAT2D} cannot hold for $G_{3,\epsilon}|_{B_0}$ (and also for $G_{3,\epsilon}|_{B_2}$ by symmetry). For $d\geq 4$, no evidence has been found of long-term trajectories contained in $\overline{A_k\cup B_k}$ for some $k\in [0,d]$. Therefore, to provide a simple systematic mechanism for the emergence of AsIUP in $G_{d,\epsilon}$ for an arbitrary $d\geq 3$, remains an inaccessible objective. 

Consequently, we target a more modest goal, which is to provide proved examples of families of piecewise $C^\infty$ symmetric maps of $S_d$, $d\geq 3$ arbitrary, that exhibit multiple acim with disjoint (asymmetric) supports when the expanding rate of their linear restrictions is sufficiently close to 1. This symmetry-breaking induced loss of ergodicity will be a consequence of the emergence of simply connected AsIUP that are generated by a multi-dimensional extension of the mechanism in $G_{2,\epsilon}$. The outcome of the construction is given in the following statement, which can be considered as the main result of this paper.
\begin{Thm}
For every $d\geq 3$, there exists $a_d\in (1,2)$ and a family $\{H_{d,a}\}_{a\in (1,,a_d)}$ of piecewise $C^\infty$ maps of $S_d$ with the following properties.
\begin{itemize}
\item The maps $H_{d,a}$ commute with the inversion symmetry $\sigma_d$ (defined in \eqref{CHANGVAR}). 
\item The maps $H_{d,a}$ coincide with $G_{d,\epsilon_a}$ on $\bigcup_{k=0}^d A_k$, where $\epsilon_a=1-\frac{a}2\in (0,\tfrac12)$ (so that the corresponding expanding rate is $a$).
\item The maps $H_{d,a}$ have an acim whose support is included in $\overline{A_k\cup B_k}$, for some $k\in \left[0,\lceil\tfrac{d}2\rceil -1\right]$. (NB: By the symmetry, they also have a disjoint acim whose support in included in $\overline{A_{d-k}\cup B_{d-k}}$). 
\item The acim support is actually included in some AsIUP, say $U$, included in $\overline{A_k\cup B_k}$. Moreover, the restriction $H_{d,a}|_{U\cap B_k}$ is an affine map. Its fixed point belongs to $B_k$ and there exists a (not necessarily normed nor orthogonal) basis $\{e_n\}_{n=1}^d$ of $\R^d$ such that the action of the linear part $L$ of $H_{d,a}|_{U\cap B_k}$ in this basis is given by 
\[
Le_n=\tfrac{a|e_n|}{|e_{n+1}|}e_{n+1},\ \forall n\in [1,d],
\]
provided that $d+1$ is identified with 1.
\end{itemize}
\label{MAINRESULT}
\end{Thm}
\begin{Rem}
\begin{itemize}
\item[(i)] The Theorem actually holds for every $k\in \left[0,\lceil\tfrac{d}2\rceil -1\right]$, with $a_d$ depending on $k$. Hence, for $a\in (1,\min_k a_d)$, it holds for all $k$. 
\item[(ii)] We do not know whether or not  the whole restriction $H_{d,a}|_{B_k}$ (and not only $H_{d,a}|_{U\cap B_k}$) can be chosen to be affine.
\end{itemize}
\end{Rem}
\noindent
{\sl Proof of Theorem \ref{MAINRESULT}:} The proof consists in providing a suitable definition of $H_{d,a}$ in $B_k$ so that we get an IUP inside $\overline{A_k\cup B_k}$. To that goal, the first step identifies an adequate collection of linearly independent points in $B_k$.
\begin{Lem}
For every $k\in [0,d]$, there exists a collection $\{p_n\}_{n=0}^d$ of linearly independent points which satisfy the following conditions
\begin{itemize}
\item $p_0$ lies in the interior of $B_k$.
\item the points $\{p_n\}_{n=1}^d$ lie in the facet common to $A_k$ and $B_k$. Moreover, $p_1$ is also included in the segment $[p_0v_k]$.
\item the vector (Euclidean) lengths $\{|p_0p_n|\}_{n=1}^{d}$ and $|p_0v_k|$ satisfy the following inequalities
\[
|p_0p_1|<|p_0p_2|<\cdots <|p_0p_d|<|p_0v_k|.
\]
\end{itemize}
\end{Lem}
\noindent
{\sl Proof.} Let $T=\overline{A_k}\cap\overline{B_k}$ be the common facet to $A_k$ and $B_k$ and let $p_1\in \text{Int} (T)$ be arbitrary. Let $p_0\in (v_kp_1)\cap \text{Int}(B_k)$ where $(v_kp_1)$ is the line through $v_k$ and $p_1$. Since $p_1\in \text{Int}(T)$, we must have 
\[
|p_0p_1|<\max_{p\in T}|p_0p|,
\]
hence, by continuity there exists $p_d\in \text{Int}(T)$ such that $|p_0p_1|<|p_0p_d|$. For the remaining points $\{p_n\}_{n=2}^{d-1}$ ($d\geq 3$), we use an iterative argument. By continuity again, there exists $p_{d-1}\in \text{Int}(T)\setminus (p_1p_d)$ such that
\[
|p_0p_1|<|p_0p_{d-1}|<|p_0p_d|.
\]
In order to obtain $p_{d-1}$, it suffices to pick up a point in $\text{Int}[p_1p_d]$ such that these inequalities hold, and then to apply a small perturbation transverse to this segment. 

If $d\geq 4$, let $(p_1p_{d-1}p_d)$ denotes the plane defined by these points and let $p_{d-2}\in \text{Int}(T)\setminus (p_1p_{d-1}p_d)$ be such that
\[
|p_0p_1|<|p_0p_{d-2}|<|p_0p_{d-1}|.
\] 
As before, in order to obtain $p_{d-2}$, it suffices to pick up a point in $\text{Int}[p_1p_{d-1}]$ such that these inequalities hold, and then to apply a small perturbation transverse to the plane $(p_1p_{d-1}p_d)$. We continue likewise for the remaining points $\{p_n\}_{n=2}^{d-3}$. Notice that for the last point $p_2$, the complement space to the hyperplane $(p_1p_3\cdots p_{d-1}p_d)$ is one-dimensional. 
\hfill $\Box$

Let $C_k$ be the $d$-simplex included in $B_k$ and defined by the vertices $\{p_n\}_{n=0}^d$. In order to construct an IUP of $H_{d,a}$ that contains an acim, it suffices to specify the (linear) action of $H_{d,a}|_{B_k}$ on $C_k$. The action on the complementary set $B_k\setminus \overline{C_k}$ is irrelevant for our purpose. It only needs to be defined such that $H_{d,a}|_{B_k} B_k\subset S_d$, in order to ensure that the whole dynamics is well-defined as a map of $S_d$ into itself.
  
\begin{Def}
Given $k\in [0,d]$, let $\{p_n\}_{n=0}^d$ be a collection as in the previous statement and let $C_k$ be the simplex defined by these points. Let $a'_d>1$ be such that 
\[
a'_d|p_0p_n|\leq |p_0p_{n+1}|\ \text{for}\ n\in [1,d-1]\quad \text{and}\quad a'_d|p_0p_d|\leq |p_0v_k|.
\]
Let $\{H_{d,a}\}_{a\in (1,a'_d)}$ be a family of piecewise affine maps defined on $A_k\cup C_k$ as follows 
\begin{itemize}
\item $H_{d,a}|_{A_k}=G_{d,1-\frac{a}2}|_{A_k}$.
\item $H_{d,a}|_{C_k}$ is an affine map with fixed point $p_0$ and with 
the following features
\[
H_{d,a}|_{C_k}p_0p_n=a\tfrac{|p_0p_n|}{|p_0p_{n+1}|}p_0p_{n+1}\ \text{for}\ n\in [1,d-1]\quad \text{and}\quad H_{d,a}|_{C_k}p_0p_d=a\tfrac{|p_0p_d|}{|p_0p_{1}|}p_0p_{1}.
\]
\end{itemize}
\label{DEFHd}
\end{Def}
Notice that the restriction $H_{d,a}|_{C_k}$ has similar permutation-expanding features as the map $G_{2,1-\frac{a}2}|_B$ in Section \ref{S-G2}. Only the length ratios $\tfrac{|H_{d,a}p_0p_n|}{|p_0p_n|}$ differ from those of $\tfrac{|G_{2,1-\frac{a}2}p_0p_n|}{|p_0p_n|}$. Yet, these choice do not really matter for our purpose as long as the product over a cycle is larger than one (so that the iterate $(H_{d,a}|_{C_k})^d$ is expanding - see the end of the proof below) and $H_{d,a}C_k\in S_d$. The features in Definition \ref{DEFHd} imply the existence of an IUP for $H_{d,a}$, as claimed in the following statement. 
\begin{Pro}
Let $\{H_{d,a}\}_{a\in (1,a'_d)}$ be as in Definition \ref{DEFHd}. Then
\begin{itemize}
\item[(i)] for all $a\in (1,a'_d)$, we have $H_{d,a}C_k\subset \overline{A_k\cup C_k}$,
\item[(ii)] when $a$ is sufficiently close to 1, the set $C_k\cup H_{d,a}C_k$ is an IUP of $H_{d,a}$. More precisely, we have $H_{d,a}(H_{d,a}C_k\cap A_k)\subset H_{d,a}C_k$.
\end{itemize}
\label{PROHd}
\end{Pro}
\noindent
{\sl Proof.} {\em (i)} That $H_{d,a}|_{C_k}$ is affine and non-singular implies that $H_{d,a}C_k$ is a $d$-simplex. Moreover, the conditions in the definition of $H_{d,a}|_{C_k}$ imply that its vertices must satisfy the conditions
\begin{equation}
H_{d,a}|_{C_k}p_n\in B_k,\ \forall n\in [0,d-1]\quad \text{and}\quad H_{d,a}|_{C_k}p_d\in A_k.
\label{LOCATION}
\end{equation}
By convexity of $\overline{A_k\cup B_k}$, it follows that $H_{d,a}C_k\subset \overline{A_k\cup B_k}$. Moreover, we have
\[
H_{d,a}|_{C_k}p_n\in [p_0p_{n+1}],\ \forall n\in [0,d-1]\quad\text{and}\quad H_{d,a}|_{C_k}[p_0p_d]\cap B_k=[p_0p_1],
\]
which, by convexity, implies $H_{d,a}C_k\cap B_k\subset C_k$. Statement {\em (i)} immediately follows.

\noindent
{\em (ii)} We begin with the following assertion.
\begin{Claim}
The set $H_{d,a}C_k\cap A_k$ is a $d$-simplex.
\end{Claim}
\noindent
{\sl Proof of the Claim.} According to \eqref{LOCATION}, each segment $H_{d,a}|_{C_k}[p_dp_n]$ $(n\in [0,d-1])$ intersects the facet $\overline{A_k}\cap \overline{B_k}$. Therefore, the set $H_{d,a}C_k\cap A_k$ can be regarded as the truncation of the polyhedral sector defined by the rays $H_{d,a}|_{C_k}[p_dp_n)$ $(n\in [0,d-1])$ by the hyperplane associated with $\overline{A_k}\cap \overline{B_k}$. As such, it must be a $d$-simplex. \hfill $\Box$

\begin{Claim}
When $a$ is sufficiently close to 1, we have $H_{d,a}(H_{d,a}C_k\cap A_k)\subset H_{d,a}C_k$.
\end{Claim}
\noindent
{\sl Proof of the Claim.} By the previous claim and the fact that $H_{d,a}|_{A_k}$ is affine and expanding, the set $H_{d,a}(H_{d,a}C_k\cap A_k)$ must be a $d$-simplex. In order to prove the claim, we study the location of the images under $H_{d,a}|_{A_k}$ of the vertices of $H_{d,a}C_k\cap A_k$, namely of the points $p_1$, $\{q_n\}_{n=1}^{d-1}$ where $q_n$ is the intersection point of the edge $H_{d,a}|_{C_k}[p_dp_n]$ and the facet $\overline{A_k}\cap \overline{B_k}$,\footnote{Notice that $p_1$ is the intersection point of the edge $H_{d,a}|_{C_k}[p_dp_0]$ and the facet $\overline{A_k}\cap \overline{B_k}$.} and $H_{d,a}|_{C_k}p_d$. 
 
By the definition of $G_{d,1-\frac{a}2}|_{A_k}$, both $H_{d,a}|_{A_k}p_1$ and $H_{d,a}|_{A_k}\circ H_{d,a}|_{C_k}p_d$ belong to the ray $[H_{d,a}|_{C_k}p_d,p_0)$. Moreover, by continuity,  when $a$ is close enough to 1, they must belong to the segment $[H_{d,a}|_{C_k}p_d,p_0]$, which is an edge of $H_{d,a}C_k$.

\begin{figure}[ht]
\begin{center}
\includegraphics*[width=38mm]{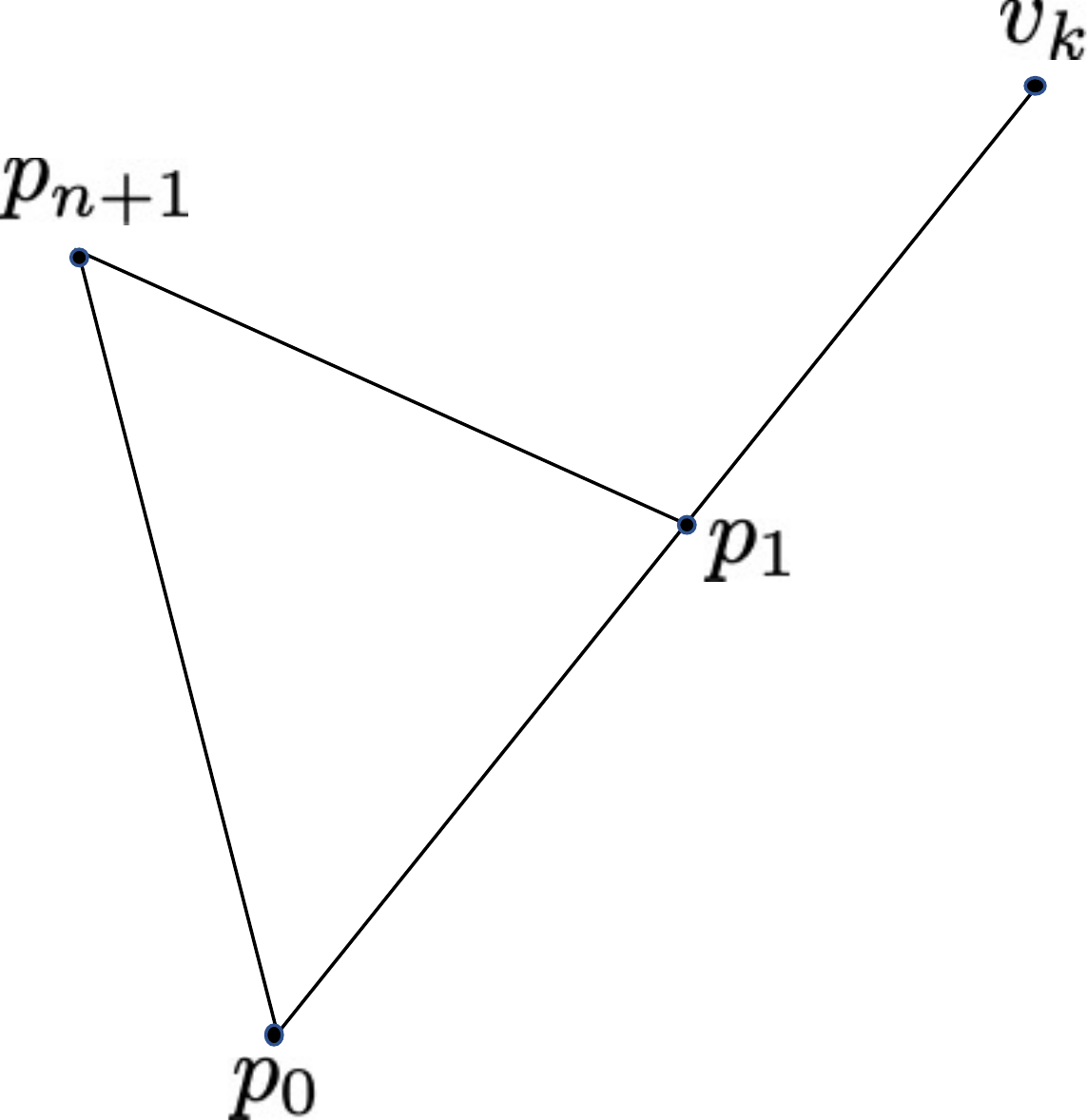}
\hspace{0.1cm}
\includegraphics*[width=53mm]{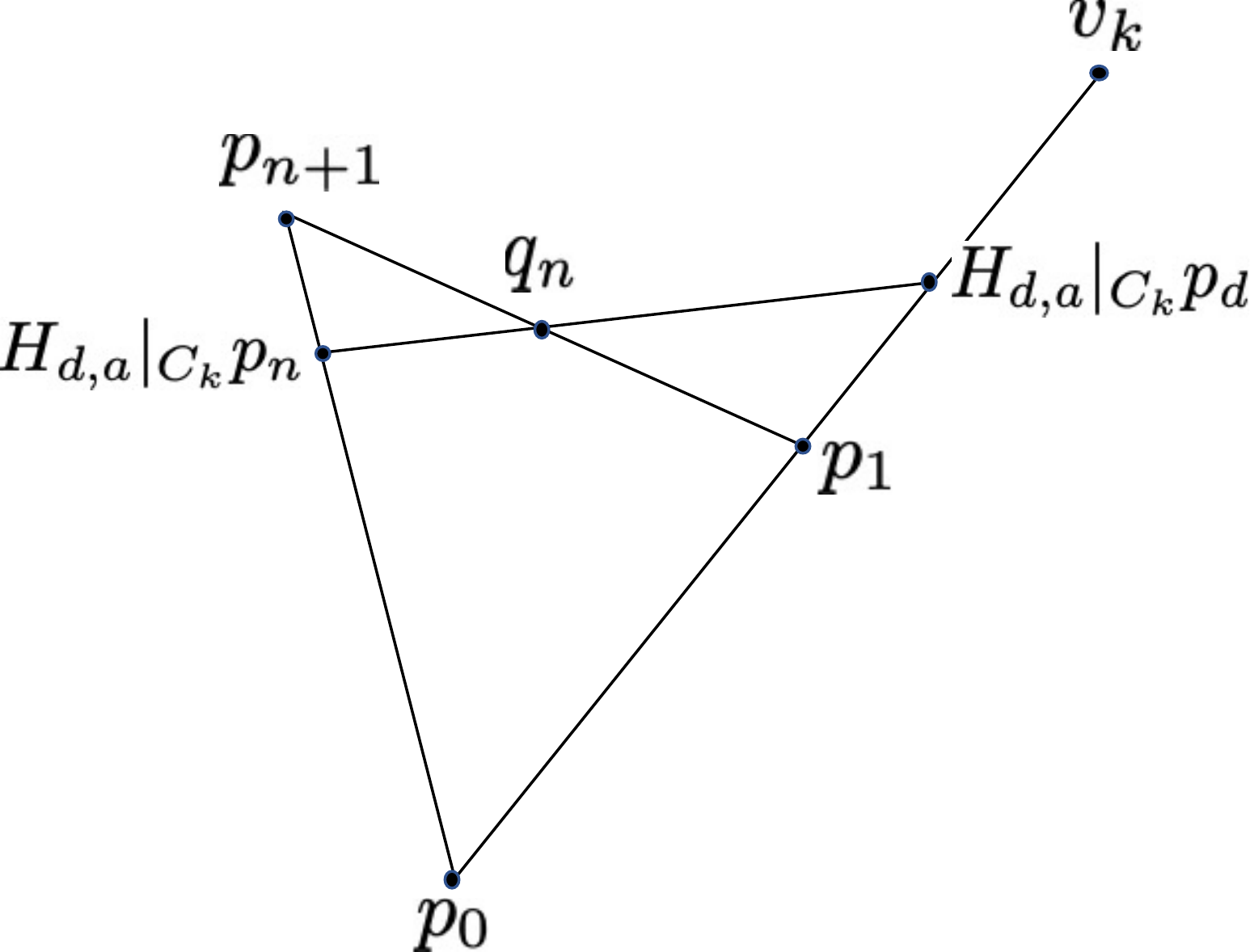}
\hspace{0.1cm}
\includegraphics*[width=63mm]{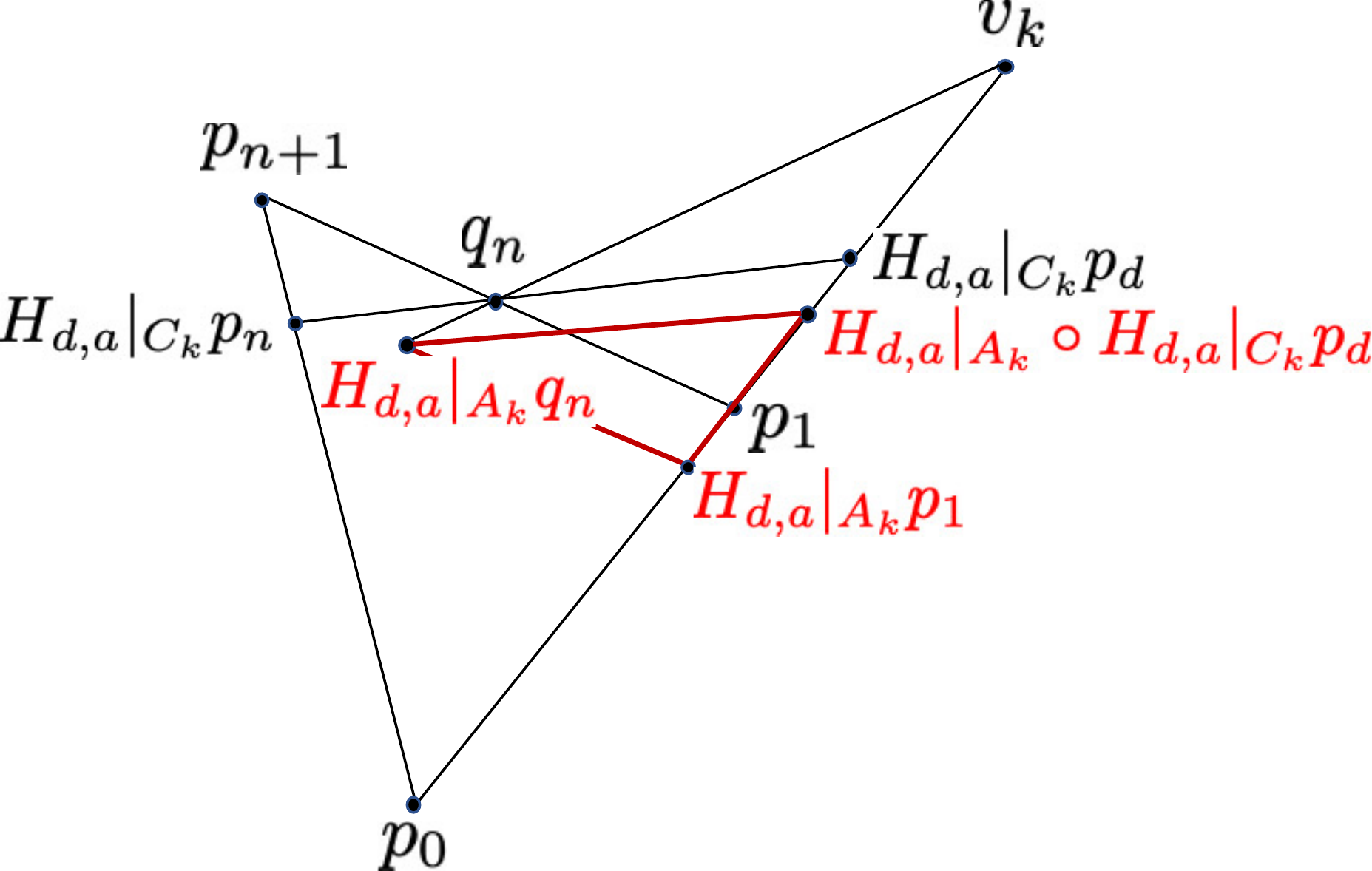}
\end{center}
\caption{Illustration of the vertices of the facets of $C_k$, $H_{d,a}C_k$ and $H_{d,a}(H_{d,a}C_k\cap A_k)$ in the 2-dimensional plane generated by the lines $(p_0v_k)$ and $(p_0p_{n+1})$. {\sl Left:} The facet of $C_k$ is the triangle $(p_0p_1p_{n+1})$. {\sl Center:}\ The facet of $H_{d,a}C_k$ is the triangle $(p_0,H_{d,a}|_{C_k}p_d,H_{d,a}|_{C_k}p_n)$. {\sl Right:}\ The facet of $H_{d,a}(H_{d,a}C_k\cap A_k)$ is the triangle $(H_{d,a}|_{A_k}p_1,H_{d,a}|_{A_k}\circ H_{d,a}|_{C_k}p_d,H_{d,a}|_{A_k}q_n)$.}
\label{DYNAPLANE}
\end{figure}
In order to locate the remaining vertices $\{H_{d,a}|_{A_k}q_n\}_{n=1}^{d-1}$, for each $n\in [1,d-1]$, we consider the 2-dimensional plane generated by the lines $(p_0v_k)$ and $(p_0p_{n+1})$, see Fig. \ref{DYNAPLANE}. Clearly, the point $q_n$ lies at the intersection of the segments $[p_1p_{n+1}]$ and $H_{d,a}|_{C_k}[p_dp_n]$. In particular, this point belongs to the edge $H_{d,a}|_{C_k}[p_dp_n]$ - but it is not a vertex - of the triangle $(p_0,H_{d,a}|_{C_k}p_d,H_{d,a}|_{C_k}p_n)$, which is a 2-facet of $H_{d,a}C_k$. Now, the definition of $H_{d,a}|_{A_k}$ and Lemma \ref{FEATURESCOUPLEDMAPS} imply that $H_{d,a}|_{A_k}q_n$ belongs to the same plane, and more precisely, to the half-plane delimited by the line $H_{d,a}|_{C_k}(p_dp_n)$ and that contains $p_0$. Therefore, by continuity, when $a$ is sufficiently close to 1, $H_{d,a}|_{A_k}q_n$ must belong to the interior of the triangle $(p_0,H_{d,a}|_{C_k}p_d,H_{d,a}|_{C_k}p_n)$.

Let finally $a$ be close enough to 1, so that all constraints above simultaneously hold. Then all points $H_{d,a}|_{A_k}p_1$, $\{H_{d,a}|_{A_k}q_n\}_{n=1}^{d-1}$ and $H_{d,a}|_{A_k}\circ H_{d,a}|_{C_k}p_d$ belong to $\overline{H_{d,a}C_k}$. The claim then follows by convexity. \hfill $\Box$

The proof of Proposition \ref{PROHd} is complete. \hfill $\Box$

Finally, that the linear part $L$ of $H_{d,a}|_{C_k}$ is as claimed in Theorem \ref{MAINRESULT} immediately follows for the basis $\{p_0p_n\}_{,=1}^d$ from the Definition \ref{DEFHd}. As a consequence, the iterated map $(H_{d,a}|_{C_k})^d$ has linear part $L^d=a^d\text{Id}$ and thus is expanding. Using also that $H_{d,a}|_{A_k}$ is expanding, we conclude that $H_{d,a}$ must have an acim supported in $C_k\cup H_{d,a}C_k$. The proof of Theorem \ref{MAINRESULT} is complete. \hfill $\Box$

\medskip

\noindent
{\bf Acknowledgements}

\noindent
We are grateful to Ilia Smilga for stimulating discussions in the early phase of this project, and to No\'e Cuneo and Matteo Tanzi for critical readings of the manuscript.

\appendix

\section{Essentials of expanding piecewise affine maps}\label{A-EPWAMAPS}
Given $d\in\N$, let $M\subset \R^d$ be a bounded polytope or $M=\T^d$. A map $F$ is said to be a {\bf piecewise affine map} of $M$ if there exists a finite collection $\{A_\omega\}$ of open, convex and disjoint polytopes included in $M$ (and called {\bf atoms}) such that
\begin{itemize}
\item the difference set $M\setminus\bigcup_{\omega}A_\omega$ has zero Lebesgue measure,
\item for each $\omega$, the restriction $F|_{A_\omega}$ is an affine map and $F(A_\omega)\subset M$. 
\end{itemize}
Notice that the {\bf atomic collection} $\{A_\omega\}$ is not unique. One may choose any finite sub-collection $\{A'_{\omega'}\}$ of atoms such that
\begin{equation}
\bigcup_{\omega'}A'_{\omega'}\subset \bigcup_\omega A_\omega\quad\text{and}\quad \text{Leb}\left(\bigcup_\omega A_\omega\setminus \bigcup_{\omega'}A'_{\omega'}\right)=0.
\label{ATOMCOLL}
\end{equation}

A piecewise affine map is said to be {\bf expanding} if there exists $a>1$ such that the linear maps $L_\omega$ associated with the affine restrictions $F|_{A_\omega}$ all satisfy the following inequality on the Euclidean lengths
\[
|L_\omega x|\geq a|x|,\ \forall x\in\R^d.
\] 
The largest of such $a$ is called the {\bf expanding rate}.

A finite union $\bigcup_{k}U_k$ of polytopes in $M$ is called an {\bf invariant union of polytopes} (IUP) for $F$ if there is an atomic collection $\{A_\omega\}$ such that
\[
F\left(\bigcup_{k,\omega}U_k\cap A_\omega \right)\subset \bigcup_{k}U_k.
\]
If $\bigcup_{k}U_k$ is an IUP of an expanding piecewise affine $F$, then $F$ must have an acim with support included in $\bigcup_{k}U_k$ \cite{T01}. Notice that this property and the acim are independent of the action of $F$ on $\bigcup_{k,\omega}U_k\setminus A_\omega$, nor they depend on the choice of the atomic collection as in \eqref{ATOMCOLL}.

An invertible transformation $\sigma:M\circlearrowleft$ is said to be an {\bf inversion symmetry} if we have $\sigma M=M$ and $\sigma^2=\text{Id}\neq \sigma$. An IUP $\bigcup_{k}U_k$ is said to be an {\bf asymmetric IUP} (AsIUP) if 
\[
\bigcup_{k}U_k\cap \sigma\left(\bigcup_{k}U_k\right)=\emptyset.
\]
If an expanding piecewise affine $F$ commutes with some inversion symmetry $\sigma$ and if $\bigcup_{k}U_k$ is an AsIUP of $F$, then evidently, $F$ must have two acim with disjoint supports.

\section{Main features of the maps $F_{N,\epsilon}$}\label{A-PHENO}
\subsection{Adapted representation of points in $\T^N$}
The coupled map phenomenology and its symmetry-induced loss of ergodicity can be more easily apprehended using the following decomposition of points in phase space. Given $u\in \R^N$, let 
\begin{itemize}
\item let $u_\text{Diag}\in \R^N$ be the vector whose coordinates are all equal to $\sum_{i=1}^Nu_i$ and
\item let $u_\perp=u-\tfrac1{N}u_\text{Diag}\in \left\{u\in\R^N\ :\ \sum_{i=1}^Nu_i=0\right\}$. 
\end{itemize}
This decomposition extends to the torus $\T^N$ as follows. If $\mathrm{u}\sim \mathrm{v}$ are two elements of the same equivalence class in $\T^N$, then we can write $\mathrm{u}=\mathrm{u}_\text{Diag}+\mathrm{u}_\perp$ and $\mathrm{v}=\mathrm{v}_\text{Diag}+\mathrm{v}_\perp$ where  
\begin{itemize}
\item $\mathrm{u}_\text{Diag}\sim \mathrm{v}_\text{Diag}$ are two elements of the same equivalent class in $\T$ and 
\item $\mathrm{u}_\perp\sim\mathrm{v}_\perp$ where this equivalence is defined by the following relation
\[
\mathrm{u}_\perp\sim\mathrm{v}_\perp\ \Longleftrightarrow\ (\mathrm{u}_\perp)_i=(\mathrm{v}_\perp)_i+n_i-\tfrac1{N}\sum_{j=1}^Nn_j,\ \forall i\in [1,N],\ 
\]
for some $n=(n_i)_{i=1}^N\in \Z^N$. Let $\mathbb{D}_N$ be the set of all equivalent classes induced by this definition. 
\end{itemize}
If $\{\mathrm{u}^t\}_{t\in\N}$ where $\mathrm{u}^{t+1}=F_{N,\epsilon}\mathrm{u}^t$ for all $t\in\N$ is an orbit of $F_{N,\epsilon}$ issued from $\mathrm{u}\in\T^N$, then the iterates $\mathrm{u}_\text{Diag}^t\in \T$ evolve independently according to the one-dimensional map $x\mapsto 2x\ \mathrm{mod}\ 1$, for which the Lebesgue measure is ergodic. Hence, any loss-of-ergodicity feature of $F_{N,\epsilon}$ has to take place in $\mathbb{D}_N$. 

Now, every element of $\mathbb{D}_N$ can be represented by an element of the scaled-centred permutahedron $P_N$ defined by 
\[
P_N=\left\{u\in \R^N\ :\ \sum_{i=1}^Nu_i=0\ \text{and}\ \sum_{i\in S}u_i\leq \tfrac{|S|(N-|S|)}{2N},\ \forall S\subsetneq [1,N],\ S\neq \emptyset\right\}.
\]
To see this, notice that the hyperplane $\sum_{i=1}^Nu_i=\tfrac{N(N+1)}2$ can be tiled by copies of the (original) permutahedron that are generated by the translations of the vectors $n\in\Z^N$ whose coordinates $n_i$ are all equal {\sl modulo} $N$ and satisfy $\sum_{i=1}^Nn_i=0$. The definition of $P_N$ then follows from the analytic characterisation of the permutahedron in terms of inequalities constraints of the coordinates, see e.g.\ Chapter 7 in \cite{LS18}, together with the appropriate scaling $\tfrac1{N}$ of the translations in the definition above of the equivalence classes associated with $u_\perp$. 
 
Furthermore, the symmetries of $F_{N,\epsilon}$ are conveyed to $P_N$, namely the map induced by $F_{N,\epsilon}$ on $P_N$ commutes with the permutations of coordinates and their sign inversion.  

\subsection{Main features of the phenomenology of $F_{N,\epsilon}$} 
As mentioned in the introduction, the coupled map phenomenology has been largely reported previously, see in particular \cite{F14,F20,S18,SB16}. Here, we provide a brief summary report together with new illustrations using the symmetric components $u_\perp$ in the permutahedrons $P_N$. 
\begin{figure}[ht]
\begin{center}
\includegraphics*[width=50mm]{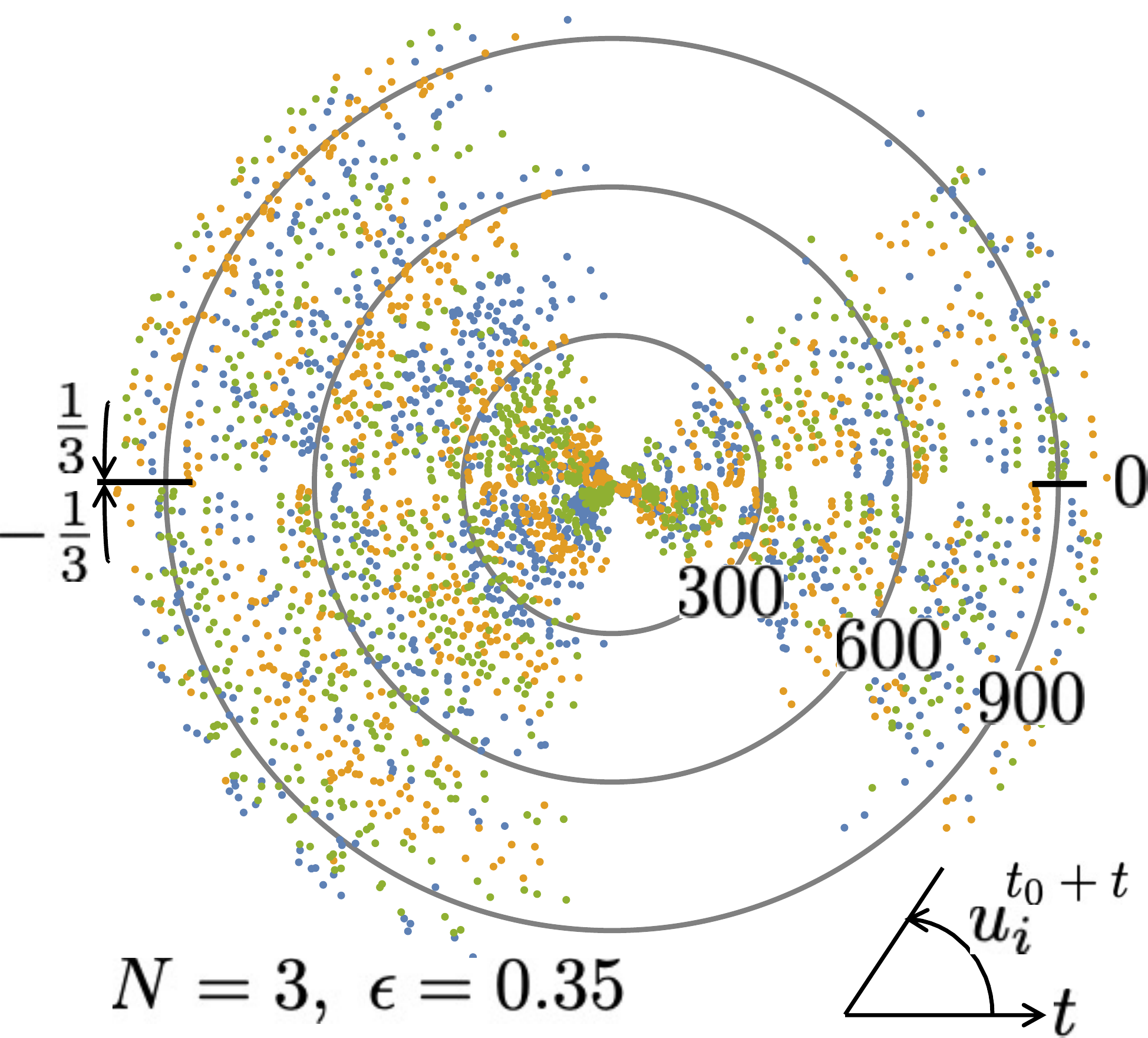}
\hspace{1cm}
\includegraphics*[width=50mm]{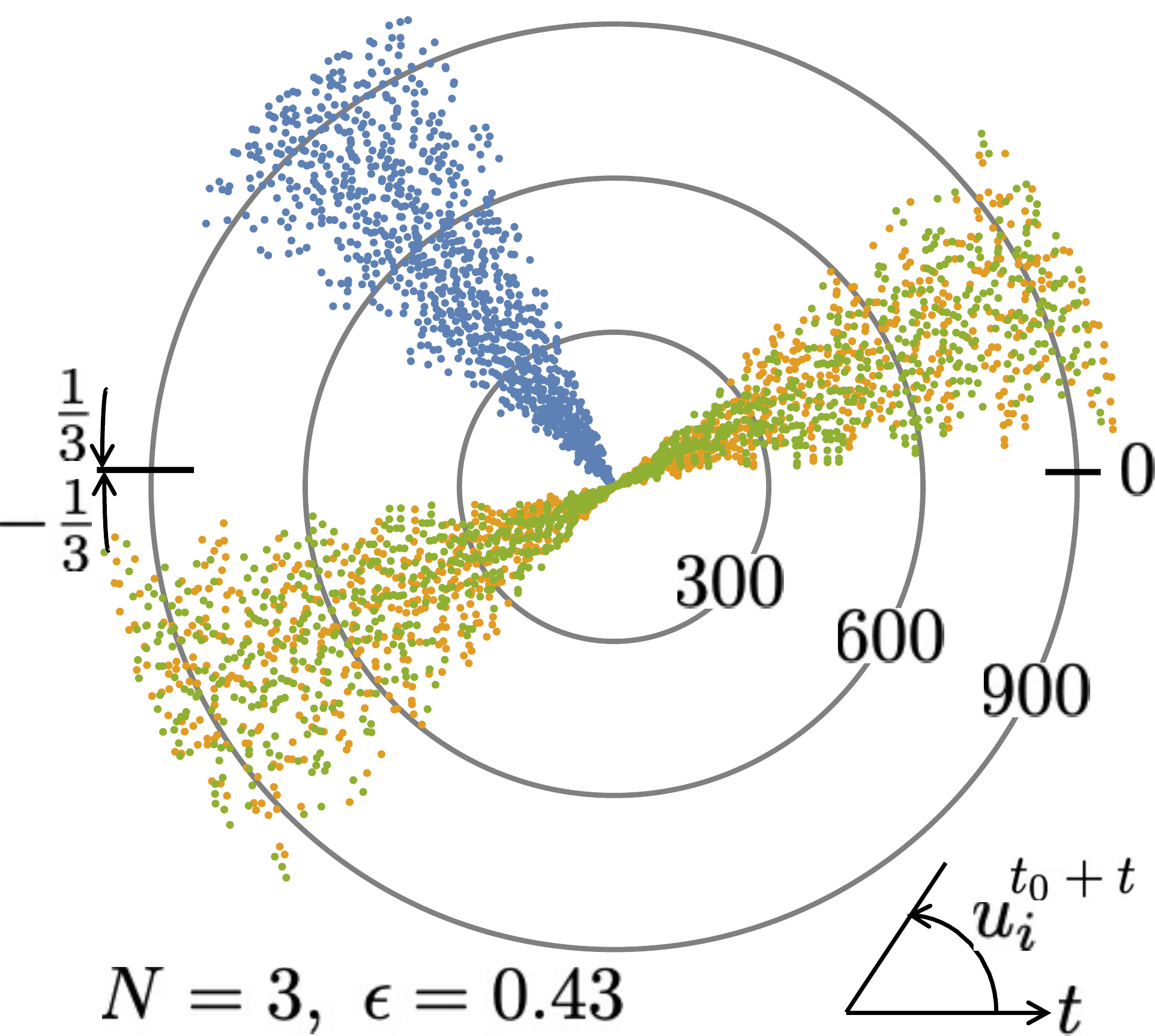}
\end{center}
\caption{Polar plots of the iterates $\{u^{t_0+t}\}_{t=1}^{1000}$ in $P_3$ ($t_0=500$), of two typical orbits of $F_{3,\epsilon}$ issued from random initial conditions $u^0$, one for $\epsilon=0.35$ in the ergodic regime (left) and one $\epsilon=0.43$ in the regime where ergodicity fails (right). Radial coordinate: $t\in [1,1000]$. Angular coordinate: $u_i^{t_0+t}$ between $-\tfrac13$ (angle = $-\pi$) and $\tfrac13$ (angle = $\pi$), one colour for each $i\in [1,3]$ (NB $\tfrac13=\tfrac{N-1}{2N}$ for $N=3$, see the definition of $P_N$ in the main text). {\sl Left.}\ Points of different colours are all scattered across various sectors, indicating that the orbit is invariant under every permutation of the three coordinates. {\sl Right.}\ The orange and green points are scattered across the same two sectors, while the blue points belong to a distinct third sector, indicating that the orbit is invariant under the permutation of two coordinates only ({\sl ie.}\ symmetry group equal to $\Pi_{i_1,i_2}$ for some pair $\{i_1,i_2\}$ of indices.). Actually, the figure suggests that the orbit generates a single connected component in phase space that also breaks the inversion symmetry $-\text{Id}|_{P_3}$, {\sl viz.}\ the reflexion with respect to the horizontal axis.}
\label{POLARPLOT1}
\end{figure}

\noindent
{\bf Notation:} Throughout this section and in Appendix \ref{A-ALTERNAT}, the symbol $\Pi_k$ ($k\in [2,N]$) denotes the {\bf group of the permutations} of the first $k$ coordinates of $u\in\R^N$ (or $\mathrm{u}\in \T^N$). The symbol $\Pi_{i_1,\cdots ,i_k}$ denotes the group of permutations of the coordinates $u_{i_1},\cdots ,u_{i_k}$ (or $\mathrm{u}_{i_1},\cdots ,\mathrm{u}_{i_k}$). 

\begin{figure}[ht]
\begin{center}
\includegraphics*[width=50mm]{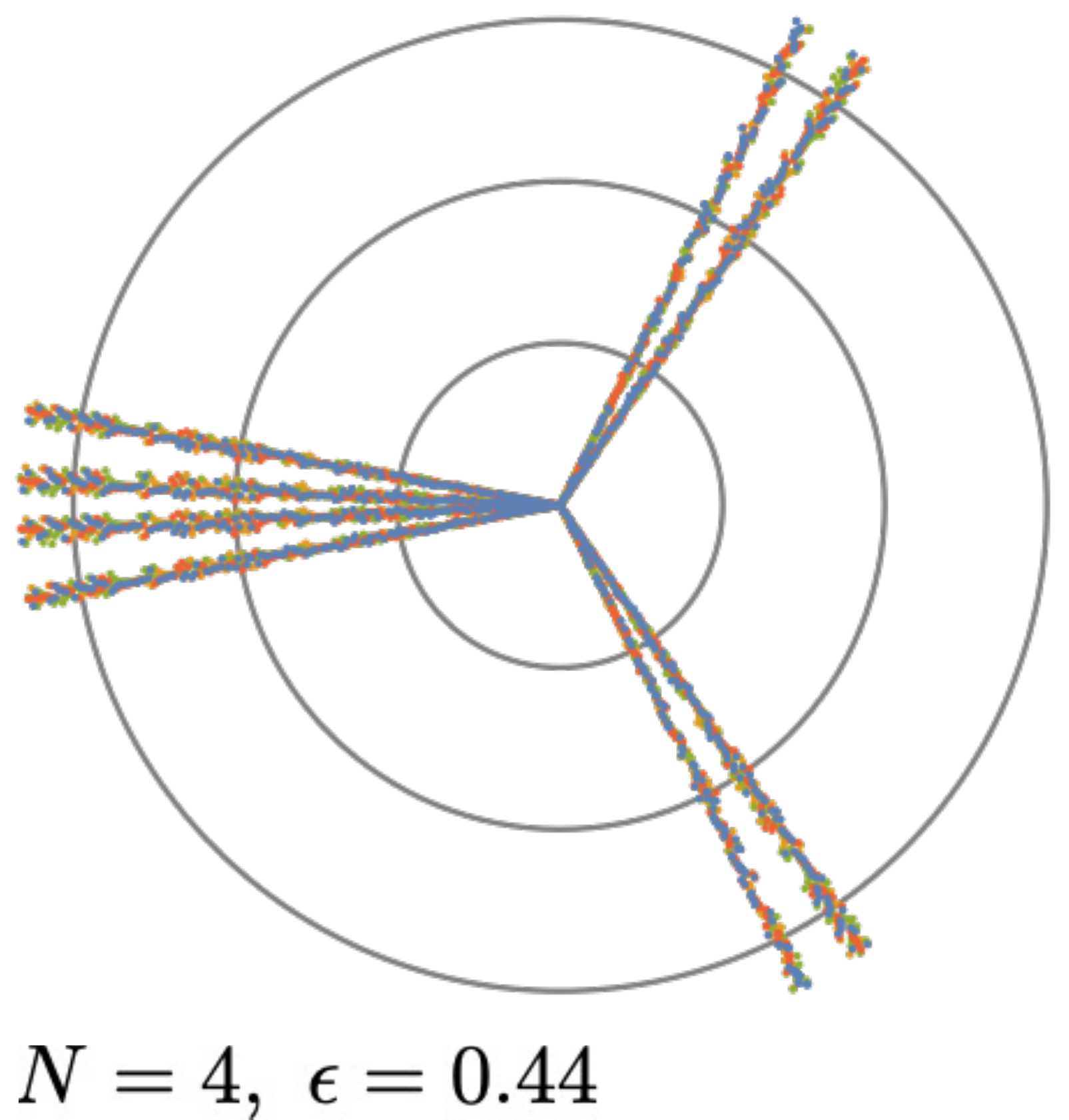}
\hspace{0.5cm}
\includegraphics*[width=50mm]{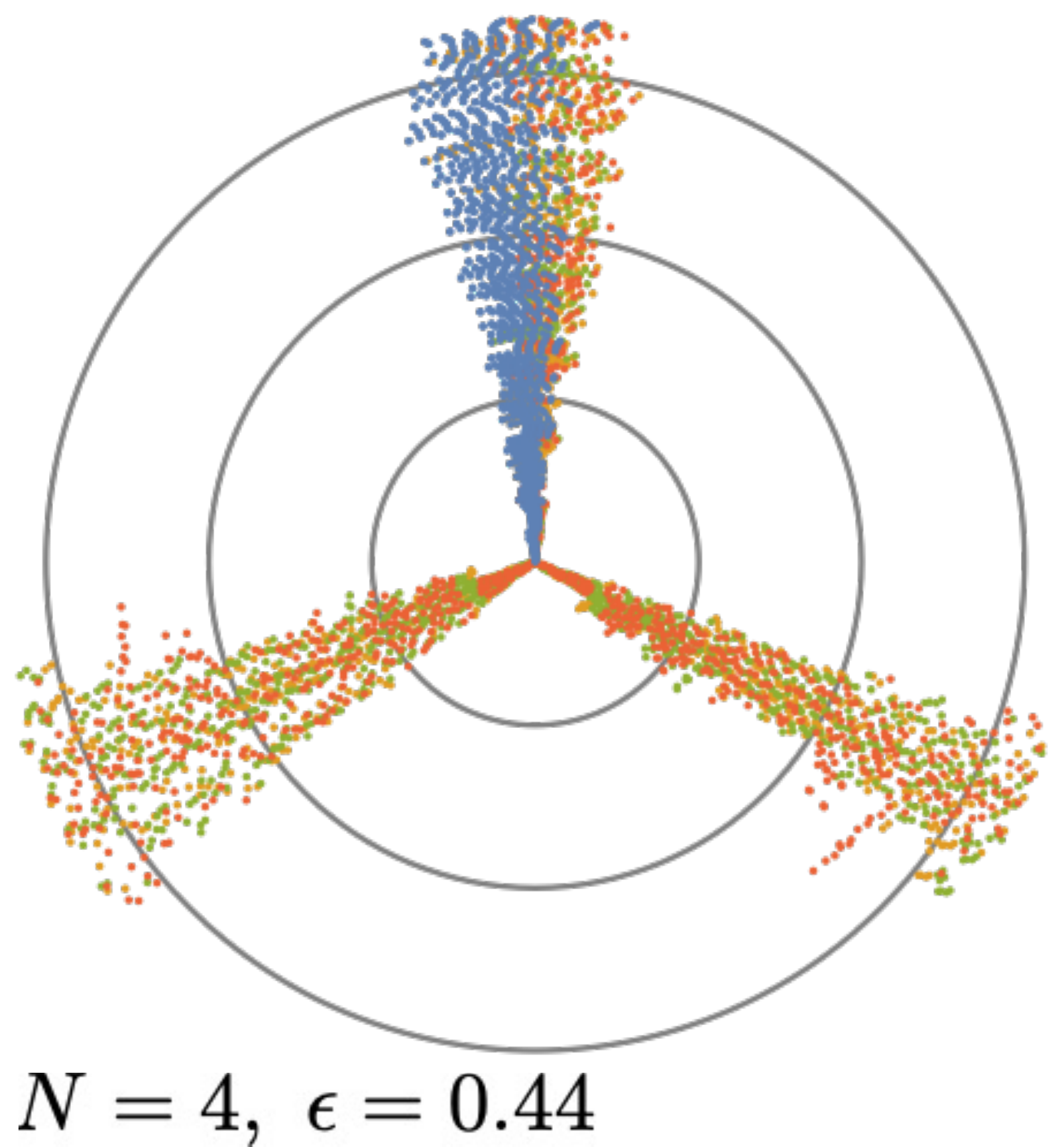}
\hspace{0.5cm}
\includegraphics*[width=50mm]{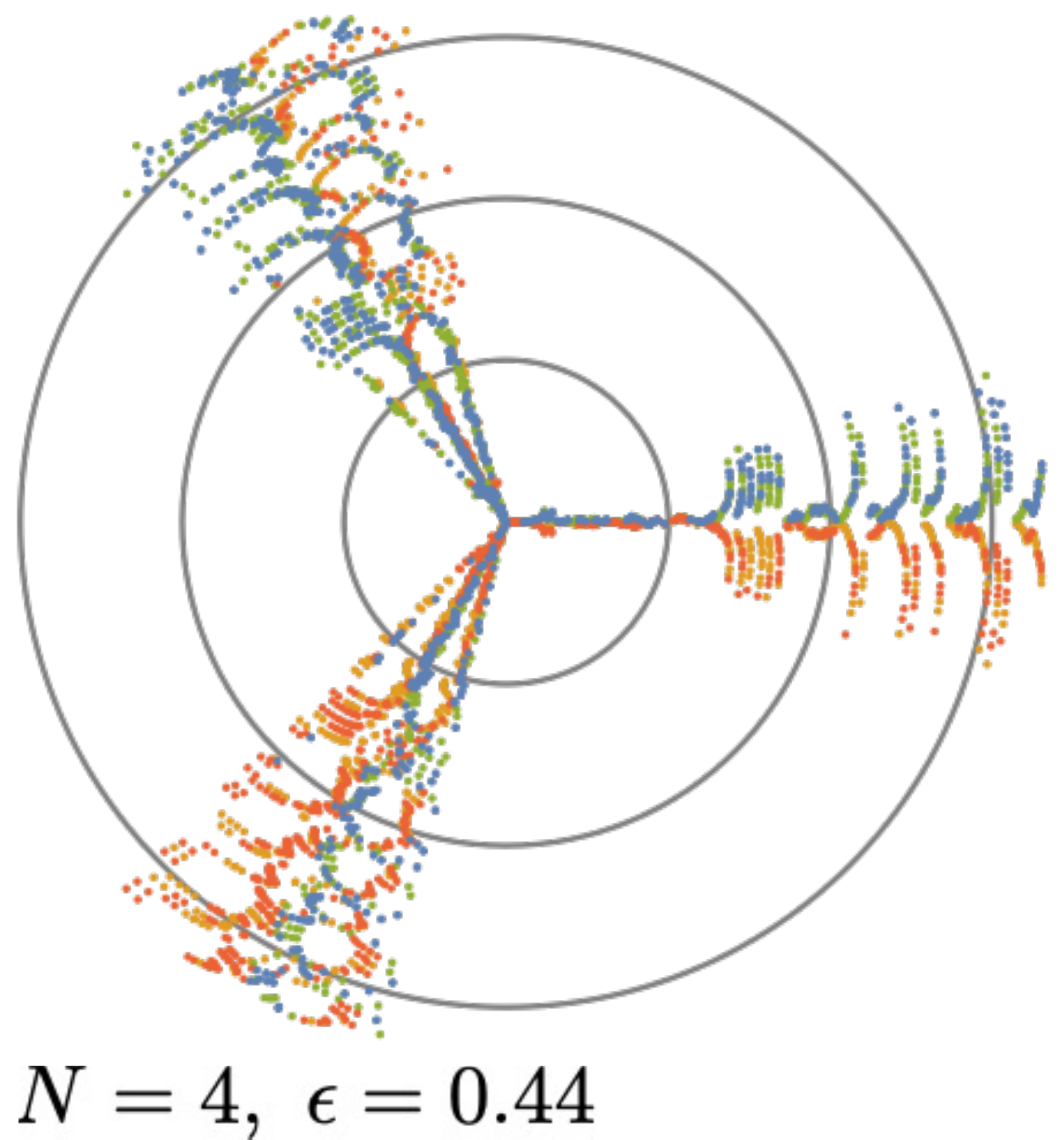}
\end{center}
\caption{Polar plots of the iterates $\{u^{t_0+t}\}_{t=1}^{1000}$ in $P_4$ of three co-existing typical orbits of $F_{4,\epsilon}$ for $\epsilon$ in the domain where ergodicity fails. Same setting as in Fig.\ \ref{POLARPLOT1}. {\sl Left.}\ Orbit invariant under every 4-coordinate permutation and $-\text{Id}|_{P_4}$ ({\sl viz.}\ symmetry group is $\Pi_4\times\Z_2$). {\sl Center.}\ Orbit invariant under every element of $\Pi_{i_1,i_2,i_3}$ for some triple $\{i_1,i_2,i_3\}$. {\sl Right.}\ Symmetry group equal to $\Pi_{i_1,i_2}\times\Pi_{i_3,i_4}$ for some permutation $\{i1,\cdots i_4\}$ of $[1,4]$.}
\label{POLARPLOT2}
\end{figure}

\begin{figure}[ht]
\begin{center}
\includegraphics*[width=50mm]{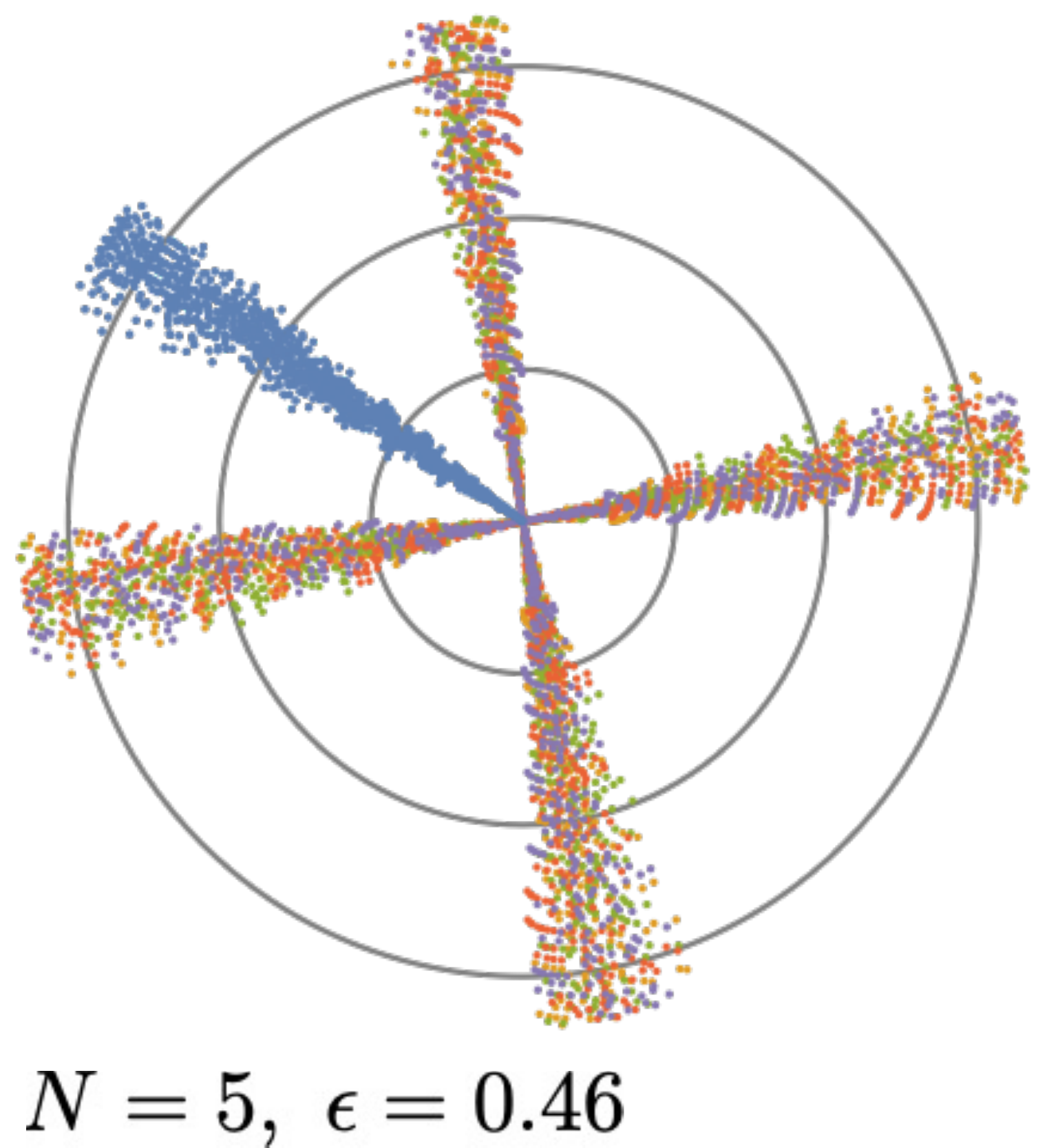}
\hspace{0.8cm}
\includegraphics*[width=50mm]{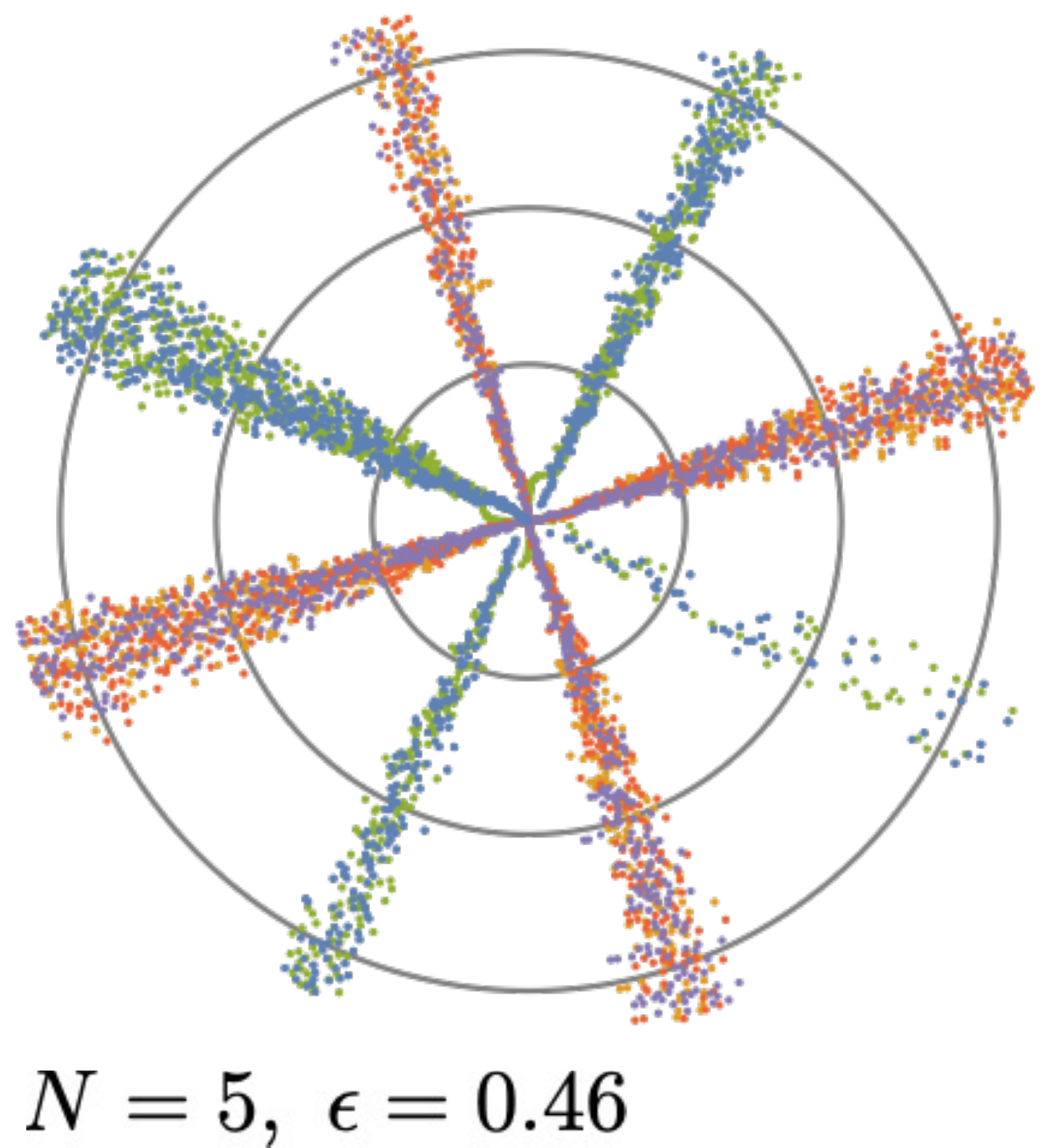}
\hspace{1cm}
\includegraphics*[width=50mm]{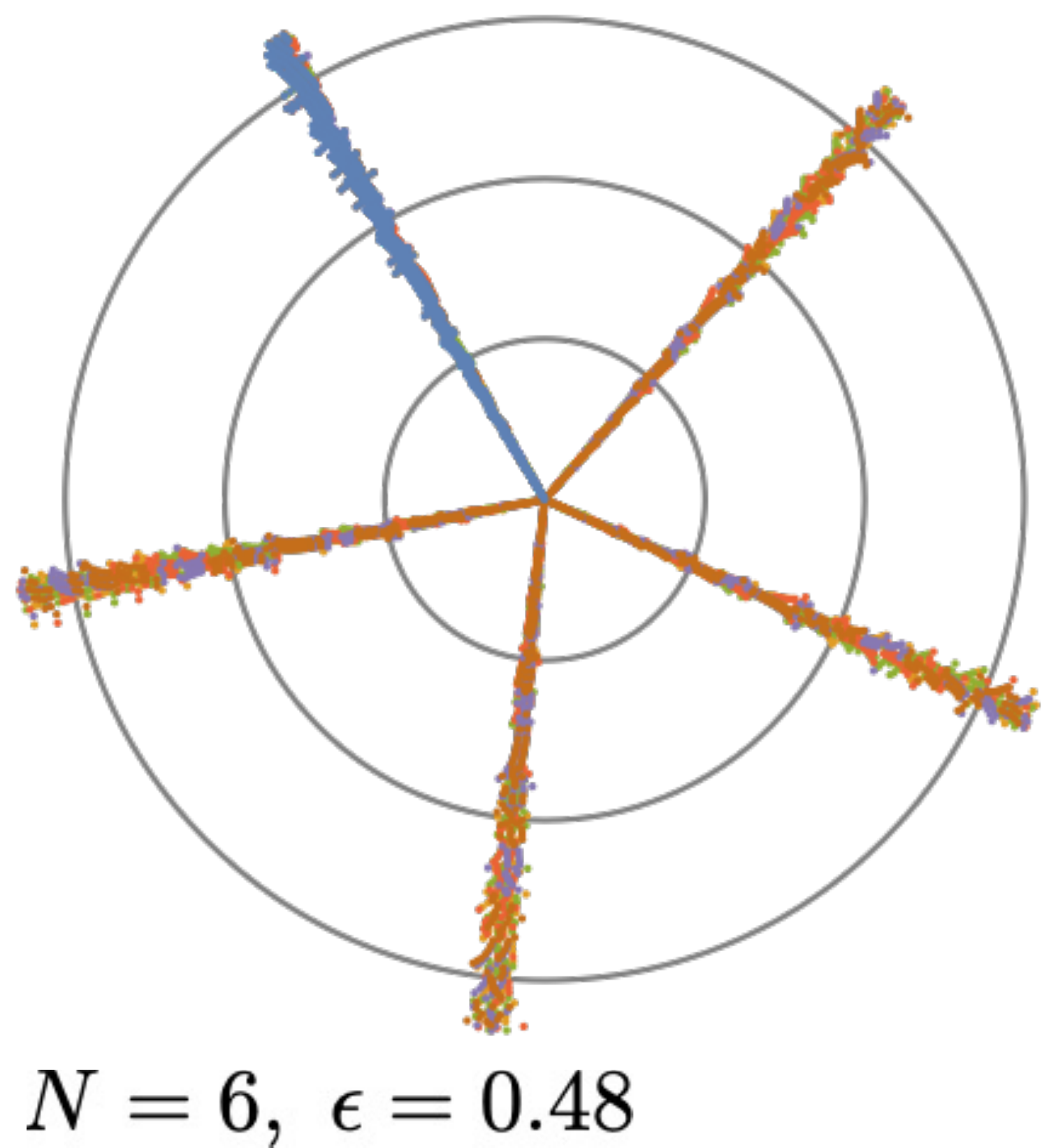}
\hspace{0.8cm}
\includegraphics*[width=50mm]{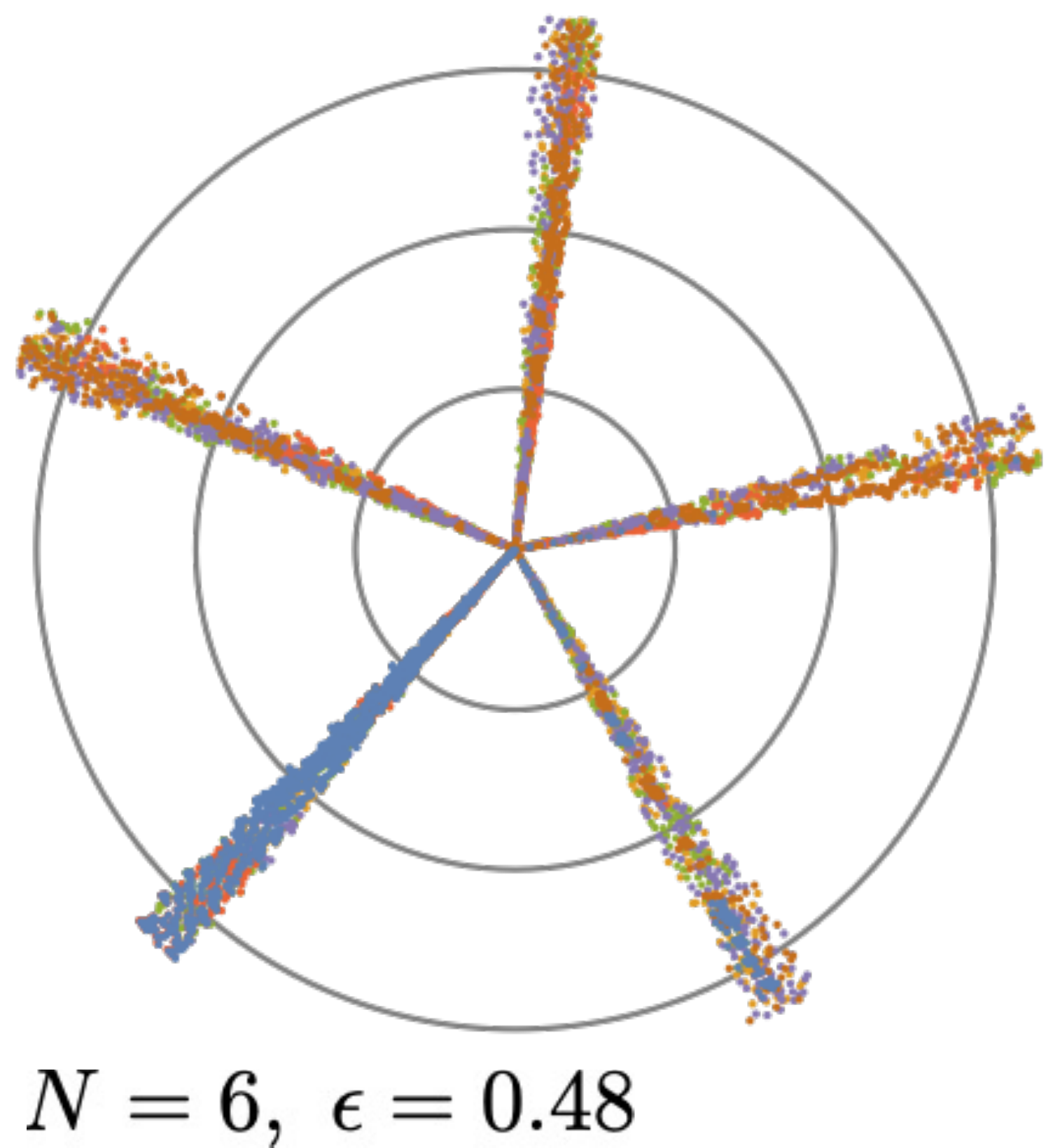}
\end{center}
\caption{Polar plots of the iterates $\{u^{t_0+t}\}_{t=1}^{1000}$ in $P_N$ of two co-existing typical orbits of $F_{N,\epsilon}$ for $\epsilon$ in the domain where ergodicity fails (top row $N=5$, bottom row $N=6$). Same setting as in Fig.\ \ref{POLARPLOT1}. {\sl Left (top and botton).}\ Orbits invariant under the elements of $\Pi_{i_1,\cdots ,i_{N-1}}$ for some $(N-1)$-uple $\{i_1,\cdots , i_{N-1}\}$. {\sl Right (top and bottom).}\ Symmetry group equal to $\Pi_{i_1,i_2}\times\Pi_{i_3,\cdots , i_N}$.}
\label{POLARPLOT3}
\end{figure}

In few words, for each $N\geq 3$,\footnote{For $N\in [1,2]$, the map $F_{N,\epsilon}$ turns out to be ergodic for all $\epsilon\in [0,\tfrac12)$.} the expanding domain $\epsilon\in [0,\tfrac12)$ can be separated into two domains, $\epsilon<\epsilon_N$ and $\epsilon>\epsilon_N$, which can be described as follows
\begin{itemize}
\item For $\epsilon<\epsilon_N$, the dynamics is ergodic, {\sl ie.}\ there is a unique ergodic component of positive Lebesgue measure. Of course, this unique component must be invariant under all symmetries, viz.\ the permutations in $\Pi_N$ and the inversion $-\textrm{Id}|_{\T^N}$. A representation of such fully symmetric component for $N=3$, obtained from numerical simulations, is given on Fig.\ \ref{POLARPLOT1} left.
\item For $\epsilon>\epsilon_N$, ergodicity is lost and there exists asymmetric (and hence multiple) ergodic components. These ergodic components have the following features
\begin{itemize}
\item For $N\geq 4$, asymmetric ergodic components may coexist with symmetric ones, see Fig.\ \ref{POLARPLOT2} for the case $N=4$. 
\item Every asymmetric component breaks the inversion symmetry $-\mathrm{Id}|_{\T^N}$, more precisely, it is disjoint from its image under $-\mathrm{Id}|_{\T^N}$. Every asymmetric component also breaks some permutation symmetry in $\Pi_N$, yet it shows a residual symmetry, {\sl ie.}\ it is invariant under some subgroup of $\Pi_N$.
\item In particular, for $\epsilon$ sufficiently close to $\tfrac12$, for some and hence every $(N-1)$-uple $\{i_1,\cdots ,i_{N-1}\}$, there exist components that are invariant under every element in $\Pi_{i_1,\cdots ,i_{N-1}}$. Examples of such components are given on Fig.\ \ref{POLARPLOT1} right ($N=3$), Fig.\ \ref{POLARPLOT2} center ($N=4$) and Fig.\ \ref{POLARPLOT3} left ($N=5,6$). 
Such ergodic components consist of $(N-1)!$ connected components, which are the images under the transformations in $\Pi_{i_1,\cdots ,i_{N-1}}$, of one of these components. This feature contrasts with those of other symmetric and asymmetric ergodic components whose connected components after identification through the action of the symmetry subgroup, do not reduce to a singleton, see Fig.\ \ref{POLARPLOT1} - \ref{POLARPLOT3} for illustrations. 
\item For $N\geq 4$, ergodic components with other residual symmetry subgroups of $\Pi_N$ such as product subgroups, may emerge (and persist) at different values of $\epsilon$. In particular, for $N=4$, ergodic components that are invariant under every transformation in $\Pi_2\times \Pi_2$ (up to conjugacy) emerge, see Fig.\ \ref{POLARPLOT2} right. For $N=5$, ergodic components with residual symmetry $\Pi_2\times \Pi_3$ emerge and for $N=6$, components with residual symmetry $\Pi_2\times \Pi_4$ has been observed, see Fig.\ \ref{POLARPLOT3} right.
\end{itemize}
\end{itemize}

\section{The dynamics of the Lorenz-type maps $f_a$}\label{A-LORENZ}
The one-dimensional maps $f_a:[0,1]\circlearrowleft$ can be characterised as follows.
\begin{itemize}
\item They commute with the reflection $x\mapsto 1-x$. 
\item They are affine maps with slope $a$ on each interval $[0,x_\mathrm{d}]$ (and then $[1-x_\mathrm{d},1]$ by symmetry) and $(x_\mathrm{d},1-x_\mathrm{d})$, where $d\in (\tfrac14,\tfrac12)$ does not depend on $a$.
\item Each of the points $0$ (and then $1$ by symmetry) and $\tfrac12$ is a fixed point of every map $f_a$.
\end{itemize}
Clearly, the condition on $d$ implies that we must have $f_a(d)\in (\tfrac12,1)$ when $a$ is sufficiently close to 2 (see Figure \ref{LORENZMAPS}). In this regime, the map $f_a$ is locally eventually onto \cite{W79} and hence ergodic with respect to some absolutely continuous measure supported on $(f_a((x_\mathrm{d})^+),f_a((1-x_\mathrm{d})^-))$. 

On the other hand, $f_a(d)\in (0,\tfrac12)$ when $a$ is sufficiently close to $1$. In this regime, both the intervals $(f_a(x_\mathrm{d}^+),f_a(x_\mathrm{d}^-))$ and $(f_a((1-x_\mathrm{d})^+),f_a((1-x_\mathrm{d})^-))$ are invariant and hence $f_a$ must have an acim on each interval (whose support being the whole interval because the map is locally eventually onto therein). Ergodicity has been lost via symmetry-breaking.

\section{Coupled map $F_{\rho,\epsilon}$ with arbitrary distribution $\rho$}\label{A-CLUSTERMAPS}
A $N$-dimensional vector $\rho=(\rho_i)_{i=1}^N$ where all $\rho_i\geq 0$ and $\sum_{i=1}^N\rho_i=1$ is called a {\bf distribution}. Given a distribution and a number $\epsilon\in [0,\tfrac12)$, consider the map $F_{\rho,\epsilon}:\T^N\circlearrowleft$ defined by \cite{FS22}
\[
(F_{\rho,\epsilon}\mathrm{u})_i=2\left(\mathrm{u}_i+\epsilon\sum_{j=1}^{N}\rho_j g(\mathrm{u}_j-\mathrm{u}_i)\right)\ \text{mod}\ 1,\ \forall i\in [1,N].
\]
All maps $F_{\rho,\epsilon}$ are expanding piecewise affine maps and their atomic partitions and symmetries are the same as those of the $F_{N,\epsilon}$. Moreover, the former are retrieved for the {\bf uniform distribution} $\rho_i=\tfrac1{N}$ for all $i$, {\sl ie.}\ we have 
\[
F_{N,\epsilon}:=F_{(\tfrac1{N})_{i=1}^N,\epsilon}.
\]
For distributions $\rho$ with rational coordinates, the maps $F_{\rho,\epsilon}$ capture the so-called cluster dynamics in the maps $F_{N,\epsilon}$, namely the dynamics in some invariant subsets of $\T^N$. To see this, given any $\mathrm{u}\in\T^N$, let the distribution $(\tfrac{n_k}{N})_{k=1}^K$ be defined by the number $K\leq N$ of groups - called clusters - inside which the coordinates $\mathrm{u}_i$ are equal, and by the number $n_k$ of coordinates in each group, see e.g.\ \cite{BKOVZ02}. The mean field coupling in $F_{N,\epsilon}$ implies that the set of configurations with given distribution $(\tfrac{n_k}{N})_{k=1}^K$ is invariant under the action of $F_{N,\epsilon}$ and the dynamics therein is governed by $F_{(\tfrac{n_k}{N})_{k=1}^K,\epsilon}$. 

Consider a distribution $\rho$ for which the first $d=N-1$ coordinates are equal, {\sl ie.}\ there exists $\varrho\in (0,\tfrac1{d})$ such that 
\[
\rho_i=\varrho,\ \forall i\in [1,d]\quad\text{and}\quad \rho_N=1-d\varrho.
\]
Evidently, the corresponding coupled maps $F_{\rho,\epsilon}$ commute with every $\pi\in\Pi_{d}$, and also with the inversion symmetry $S=-\text{Id}|_{\T^N}$. The arguments in Section \ref{S-TORUSMAPS} imply that when the corresponding projected map ${\cal F}_{\rho,\epsilon}$ has an AsIUP with respect to the inversion symmetry $\sigma_\Sigma$ defined by \eqref{INVSYM}, then $F_{\rho,\epsilon}$ must have two acim with disjoint supports.

Moreover, as in Claim \ref{CONJUG}, ${\cal F}_{\rho,\epsilon}$ is conjugated to a skew-product dynamical system whose base map, say $G_{\rho,\epsilon}$ is a mapping of $S_d$ into itself. This map commutes with $\sigma_d$ and the existence of an AsIUP for $G_{\rho,\epsilon}$ implies the existence of two acim with disjoint supports for $F_{\rho,\epsilon}$.

Now, similarly as in Lemma \ref{FEATURESCOUPLEDMAPS}, the main features of the restrictions $G_{\rho,\epsilon}|_{A_k}$ and $G_{\rho,\epsilon}|_{B_k}$ are given in the following statement.
\begin{Lem}
In addition to commuting with $\sigma_d$, the expanding piecewise affine map $G_{\rho,\epsilon}$ has the following features for every $\varrho\in \left(0,\tfrac1{d}\right)$ and $\epsilon\in \left(0,\tfrac12\right)$.

\noindent
(i) Every simplex $A_k$ is an atom of $G_{\rho,\epsilon}$. The restrictions of $G_{\rho,\epsilon}$ to $A_0$ and to $A_k$ do not depend on the value of $\varrho\in \left(0,\tfrac1{d}\right)$ and they respectively write
\[
(G_{\rho,\epsilon}|_{A_0}x)_i=2(1-\epsilon) x_i\quad\text{\em and}\quad (G_{\rho,\epsilon}|_{A_k}x)_i=2(1-\epsilon) x_i+(2\epsilon-1)\delta_{i,k},\ i\in [1,d].
\]

\noindent
(ii) Let $d\geq 2$. Given $k\in [1,d-1]$, the simplex $B_k$ is an atom of $G_{\rho,\epsilon}$ iff $d\in [2,3]$ and $\varrho\geq\tfrac14$. For $k=0$ and $k=d$, the same property holds iff $\varrho\in [\tfrac1{2d},\tfrac1{2(d-1)}]$. 
\label{A-FEATURESCOUPLEDMAPS}
\end{Lem}
 
\noindent
{\sl Proof.} Throughout the proof, we regard $F_{\rho,\epsilon}$ and $P$ as maps from $\R^N$ into itself. Let also $\mathrm{F}_{\rho,\epsilon}=P\circ F_{\rho,\epsilon}:I_N\to D_\ast^N$. A careful examination of the definition of $G_{\rho,\epsilon}$ concludes that an atom of this map is defined by the simultaneous occurrence of the following conditions for the variable $u\in I_N$ 
\begin{itemize}
\item[(a)] the collection $\left\{ \sum_{j=1}^N\lfloor u_j-u_i+\tfrac12\rfloor\right\}_{i=1}^N$ is constant,
\item[(b)] the collection $\left\{ \lfloor (F_{\rho,\epsilon} u)_N-(F_{\rho,\epsilon} u)_i\rfloor\right\}_{i=1}^{N-1}$ is constant, 
\item[(c)] the ordering of the coordinates $\{(\mathrm{F}_{\rho,\epsilon}u)_i\}_{i=1}^{N-1}$ is constant.
\end{itemize}
Clearly, these conditions only depend on the differences $\{u_{i+1}-u_i\}_{i=1}^N$, namely they are genuine conditions for the variable $x=((\phi_N u)_i)_{i=1}^{d}\in S_d$.

The rest of the proof is purely computational. The equality $u_j-u_i=\sum_{n=i}^{j-1}x_n$ implies that for $x\in A_0$, we have $\lfloor u_j-u_i+\tfrac12\rfloor =0$ for all $i,j\in [1,N]$ which immediately implies (a). Moreover, we have 
\[
(F_{\rho,\epsilon}u)_j-(F_{\rho,\epsilon}u)_i=2(1-\epsilon)(u_j-u_i),
\]
and hence 
\[
0<(F_{\rho,\epsilon}u)_j-(F_{\rho,\epsilon}u)_i<1-\epsilon, \ \forall i<j\in [1,N],
\]
which implies (property (b))
\[
\lfloor (F_{\rho,\epsilon}u)_N-(F_{\rho,\epsilon}u)_i\rfloor =0, \ \forall i\in [1,N-1]
\]
and hence, together with the expression of $(F_{\rho,\epsilon}u)_j-(F_{\rho,\epsilon}u)_i$ above, the property (c) and the expression of $G_{\rho,\epsilon}|_{A_0}$. The result for $A_d$ follow by symmetry.

Now, assume that $x\in A_k$ for some $k\in \left[1,\lceil\tfrac{d}2\rceil\right]$, the remaining cases follow by symmetry. 
\begin{figure}[ht]
\begin{center}
\includegraphics*[width=90mm]{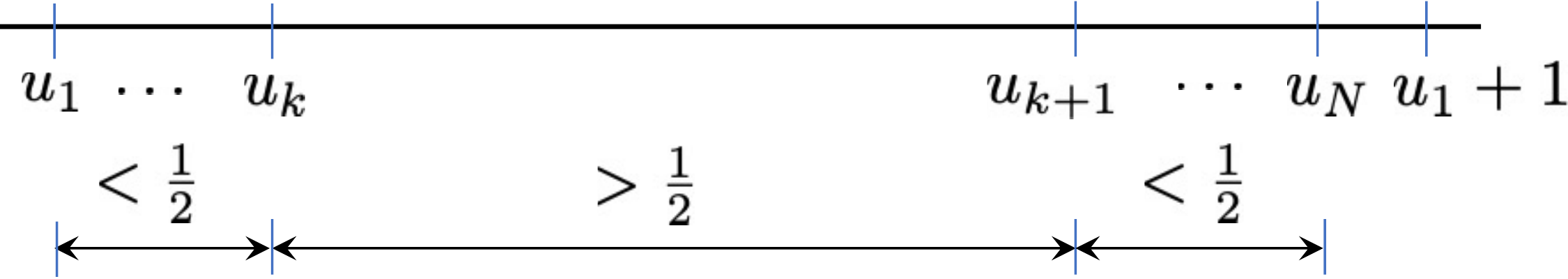}
\end{center}
\caption{Illustration of a configuration $u\in I_N$ for which $x=((\phi_N u)_i)_{i=1}^{d}\in A_k$ for $k\in [2,d]$ (to be adapted in the case $k=1$).}
\label{KI=I}
\end{figure}
Then, we have (see Fig.\ \ref{KI=I}) for $i\in [1,k]$, 
\[
\lfloor u_j-u_i+\tfrac12\rfloor=\chi_{[k+1,N]}(j),\  \forall j\in [1,N],
\]
and for $i\in [k+1,N]$
\[
\lfloor u_j-u_i+\tfrac12\rfloor=-\chi_{[1,k]}(j),\ \forall j\in [1,N].
\]
Again, property (a) is evident. Moreover, direct calculations yield
\[
(F_{\rho,\epsilon}u)_j-(F_{\rho,\epsilon}u)_i=2(1-\epsilon)(u_j-u_i)+2\epsilon \chi_{[1,k]}(i)\chi_{[k+1,N]}(j),\ \forall i<j\in [1,N],
\]
and hence (property (b))
\[
\lfloor (F_{\rho,\epsilon}u)_N-(F_{\rho,\epsilon}u)_i\rfloor =\chi_{[1,k]}(i),\ \forall i\in [1,N-1],
\]
and also the coordinates $(\mathrm{F}_{\rho,\epsilon}u)_i$ are increasing\footnote{because we have in particular $(\mathrm{F}_{\rho,\epsilon}u)_{k+1}-(\mathrm{F}_{\rho,\epsilon}u)_k>1-\epsilon+2\epsilon -1=\epsilon$.} (property (c)). The expression of $G_{\rho,\epsilon}|_{A_k}$ immediately follows.

Assume now that $x\in B_k$ for some $k\in \left[1,\lceil\tfrac{d}2\rceil\right]$.
\begin{figure}[ht]
\begin{center}
\includegraphics*[width=90mm]{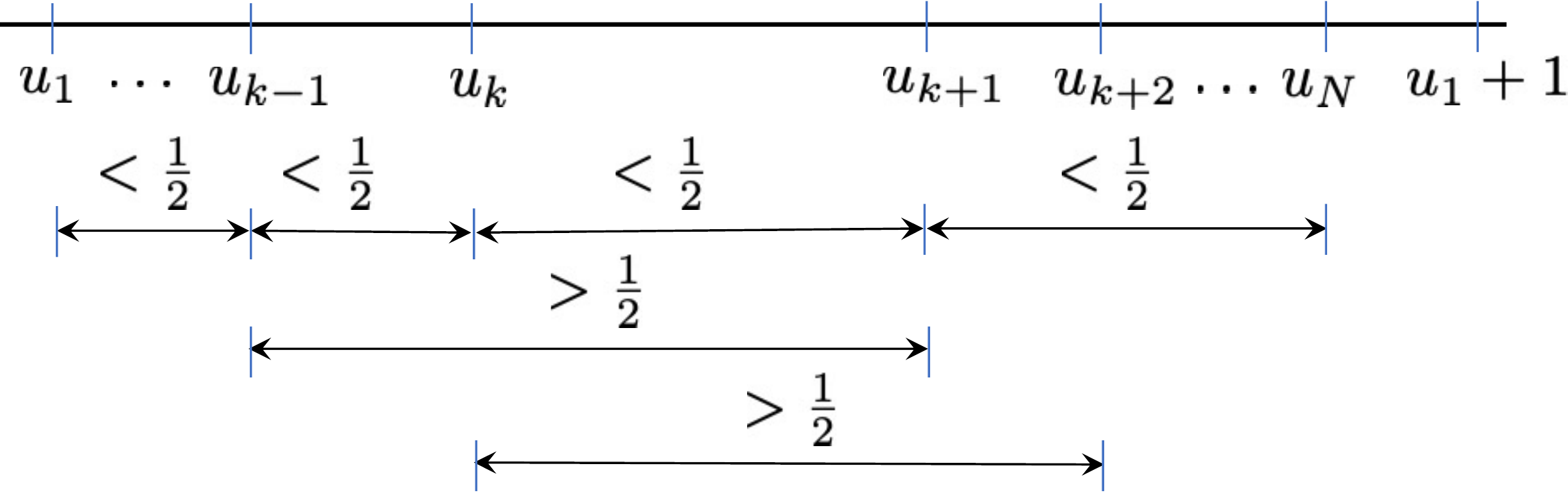}
\end{center}
\caption{Illustration of a configuration $u\in I_N$ for which $x=((\phi_N u)_i)_{i=1}^{d}\in B_k$ for $k\in [2,d-1]$ (to be adapted in the cases $k=0,1$ and $k=d$).}
\label{KI=I+1}
\end{figure}
Then, we have (see Fig.\ \ref{KI=I+1}) for $i\in [1,k-1]$\footnote{When $k=1$, we naturally ignore the indices $i\leq k-1$.}
\[
\lfloor u_j-u_i+\tfrac12\rfloor=\chi_{[k+1,N]}(j),\  \forall j\in [1,N].
\]
Moreover 
\[
\lfloor u_j-u_{k}+\tfrac12\rfloor=\chi_{[k+2,N]}(j)\ \text{and}\ \lfloor u_j-u_{k+1}+\tfrac12\rfloor=-\chi_{[1,k-1]}(j),\ \forall j\in [1,N].
\]
and for $i\in [k+2,N]$, we have
\[
\lfloor u_j-u_i+\tfrac12\rfloor=-\chi_{[1,k]}(j),\ \forall j\in [1,N].
\]
Again, property (a) is obvious. Moreover, we have
\[
\sum_{n=1}^N\rho_n\lfloor u_n-u_i+\tfrac12\rfloor-\sum_{n=1}^N\rho_n\lfloor u_n-u_N+\tfrac12\rfloor=\left\{\begin{array}{ccl}
1&\text{if}&i\in [1,k-1]\\
1-\varrho&\text{if}& i=k\\
\varrho&\text{if}&i=k+1\\
0&\text{if}&i\in [k+2,N-1]
\end{array}\right.
\]
from where we obtain (property (b))
\[
\lfloor (F_{\rho,\epsilon}u)_N-(F_{\rho,\epsilon}u)_i\rfloor =\chi_{[1,k]}(i),\ \forall i\in [1,N].
\]
As a consequence we have
\[
(\mathrm{F}_{\rho,\epsilon} u)_i=2(1-\epsilon)u_i+2\epsilon\sum_{j=1}^N\rho_j u_j-\lfloor (F_{\rho,\epsilon}u)_N\rfloor+\chi_{[1,k]}(i)+\left\{\begin{array}{ccl}
-2\epsilon\sum_{j=k+1}^N\rho_j&\text{if}&i\in [1,k-1]\\
-2\epsilon\sum_{j=k+2}^N\rho_j&\text{if}&i=k\\
2\epsilon\sum_{j=1}^{k-1}\rho_j&\text{if}&i=k+1\\
2\epsilon\sum_{j=1}^{k}\rho_j&\text{if}&i\in [k+2,N]
\end{array}\right.
\]
and hence
\[
(\mathrm{F}_{\rho,\epsilon}u)_{i+1}-(\mathrm{F}_{\rho,\epsilon}u)_i=2(1-\epsilon)(u_{i+1}-u_i)+
\left\{\begin{array}{ccl}
0&\text{if}&i\in [1,k-2]\\
2\epsilon \rho_{k+1}&\text{if}& i=k-1\\
2\epsilon (1-\rho_k-\rho_{k+1})-1&\text{if}&i=k\\
2\epsilon \rho_k&\text{if}& i=k+1\\
0&\text{if}&i\in [k+2,N-1]
\end{array}\right.
\]
It immediately follows that $(\mathrm{F}_{\rho,\epsilon} u)_{i+1}-(\mathrm{F}_{\rho,\epsilon} u)_i>0$ for all $i\neq k$. For $i=k$, the situation depends on the location of $\varrho=\rho_k=\rho_{k+1}$ with respect to $\tfrac14$, since we have using also $u_{k+1}-u_k<\tfrac12$
\[
(\mathrm{F}_{\rho,\epsilon}u)_{k+1}-(\mathrm{F}_{\rho,\epsilon}u)_k\in \left(2\epsilon (1-2\varrho)-1,\epsilon(1-4\varrho)\right),
\]
and $2\epsilon (1-2\varrho)-1\leq 0$ for all $\epsilon \in [0,\tfrac12)$ and $\varrho\leq \tfrac12$. Hence, the sign of $(\mathrm{F}_{\rho,\epsilon}u)_{k+1}-(\mathrm{F}_{\rho,\epsilon}u)_k$ is certainly negative when $\varrho\geq \tfrac14$. Otherwise, this sign depends on the location of $2(1-\epsilon)(u_{k+1}-u_k)$ with respect to $1-2\epsilon (1-2\varrho)$.

Using also that $u_{i+2}-u_i>\tfrac12$ for $i\in [k-1,k]$, we obtain
\[
(\mathrm{F}_{\rho,\epsilon}u)_{i+2}-(\mathrm{F}_{\rho,\epsilon}u)_i>\epsilon(1-2\varrho),\ i\in [k-1,k],
\]
and hence for $\varrho\geq \tfrac14$, we have (property (c))
\begin{equation}
(\mathrm{F}_{\rho,\epsilon}u)_{1}<\cdots <(\mathrm{F}_{\rho,\epsilon}u)_{k-1}<(\mathrm{F}_{\rho,\epsilon}u)_{k+1}<(\mathrm{F}_{\rho,\epsilon}u)_{k}<(\mathrm{F}_{\rho,\epsilon}u)_{k+2}<\cdots <(\mathrm{F}_{\rho,\epsilon}u)_{N-1}<(\mathrm{F}_{\rho,\epsilon}u)_{N}
\label{LASTORDER}
\end{equation}
Since we must have $\varrho<\tfrac1{N-1}$, the condition $\varrho\geq \tfrac14$ can only hold for $N\in [3,4]$. The proof of the statement {\sl (ii)} for $k\in [1,d]$ is complete. The case $k=0$ can be treated by a similar analysis. Their details are left to the reader. \hfill $\Box$ 

In addition, the ordering in equation \eqref{LASTORDER} implies that for $\varrho\geq \tfrac14$, the reduced map in $B_1$ writes
\[
(G_{\rho,\epsilon}|_{B_1}x)_i=\left\{\begin{array}{ccl}
(\mathrm{F}_{\rho,\epsilon}u)_{1}-(\mathrm{F}_{\rho,\epsilon}u)_{2}=2(1-\epsilon)(u_{1}-u_{2})+1-2\epsilon(1-2\varrho)&\text{if}&i=1\\
(\mathrm{F}_{\rho,\epsilon}u)_{3}-(\mathrm{F}_{\rho,\epsilon}u)_{1}=2(1-\epsilon)(u_{3}-u_{1})+2\epsilon (1-\varrho)-1&\text{if}&i=2\\
(\mathrm{F}_{\rho,\epsilon}u)_{4}-(\mathrm{F}_{\rho,\epsilon}u)_{3}=2(1-\epsilon)(u_{4}-u_{3})&\text{if}&N=4\ \text{and}\ i=3
\end{array}\right.
\]
{\sl ie.}
\begin{equation}
(G_{\rho,\epsilon}|_{B_1}x)_i=\left\{\begin{array}{ccl}
-2(1-\epsilon)x_1+1-2\epsilon(1-2\varrho)&\text{if}&i=1\\
2(1-\epsilon)(x_{1}+x_{2})+2\epsilon (1-\varrho)-1&\text{if}&i=2\\
2(1-\epsilon)x_{3}&\text{if}&N=4\ \text{and}\ i=3
\end{array}\right.
\label{EXPRG}
\end{equation}
Besides, the computation of $G_{\rho,\epsilon}|_{B_0}$ ($\varrho\geq\tfrac14$) for $N=4$ yields the following expression
\[
(G_{\rho,\epsilon}|_{B_0}x)_i=\left\{\begin{array}{ccl}
2(1-\epsilon)x_2&\text{if}&i=1\\
-2(1-\epsilon)(x_{1}+x_{2})+1-2\epsilon (1-3\varrho)&\text{if}&i=2\\
2(1-\epsilon)(x_{1}+x_{2}+x_3)+2\epsilon (1-2\varrho)-1&\text{if}&i=3
\end{array}\right.
\]

Let $d=2$. Lemma \ref{A-FEATURESCOUPLEDMAPS} states that the 2-dimensional map $G_{\rho,\epsilon}:S_2\circlearrowleft$ associated with any 3-dimensional distribution $\rho$ of the form $\rho=(\varrho,\varrho,1-2\varrho)$ with $\varrho\geq\tfrac14$, is an expanding piecewise affine map with atomic collection $\{A_0,A_1,A_2,B\}$ where $B=B_0=B_1=B_2$. Moreover, the restriction $G_{\rho,\epsilon}|_{B}$ has similar characteristics as those of $G_{2,\epsilon}|_{B}$ in Claim \ref{FEAT2D}.
\begin{Claim}
(i) The fixed point $p_0$ of $G_{\rho,\epsilon}|_{B}$ belongs to $B$.

\noindent
(ii) Let $p_1$ be the intersection point of the segment $[p_0v_2]$ and the edge $\overline{A_2}\cap\overline{B}$ and let $p_2$ be the intersection point (which exists) of the image segment $G_{\rho,\epsilon}|_B[p_0p_1]$ and $\overline{A_2}\cap\overline{B}$. In the basis formed by the vectors $p_0p_1$ and $p_0p_2$, the linear part of $G_{\rho,\epsilon}|_{B}$ is given by the following matrix
\[
2(1-\epsilon)\left(\begin{array}{cc}
0&\tfrac1{\alpha_\epsilon}\\
\alpha_\epsilon&0\end{array}\right)
\]
for some $\alpha_\epsilon>0$.
\end{Claim}
\noindent
{\sl Proof:} {\sl (i)} The coordinates of $p_0$ are $\left(\tfrac{1-2\epsilon(1-2\varrho)}{3-2\epsilon},\tfrac{1-2\epsilon\varrho}{3-2\epsilon}\right)$, from where one checks that $p_0\in B$ for all $\epsilon\in (0,\tfrac12)$. 

\noindent 
{\sl (ii)} The linear part of $G_{\rho,\epsilon}|_{B}$ writes $2(1-\epsilon)M$ where $M$ has eigenvector $e_2$ with eigenvalue $1$, and $2e_1-e_2$ with eigenvalue $-1$ (Recall that the $e_j$ are the canonical vectors). Using that
\[
p_0p_1=x_\epsilon e_1+y_\epsilon (2e_1-e_2),
\]
for some $y_\epsilon\neq 0$, we obtain that $M$ must write as claimed in the basis formed by $p_0p_1$ and $G_{\rho,\epsilon}|_Bp_0p_1$. Finally, one checks that the ray $[p_0,G_{\rho,\epsilon}|_Bp_1)$ intersects the edge $\overline{A_2}\cap {B}$ for all $\epsilon\in (0,\tfrac12)$. \hfill $\Box$

Similarly as in the case $\varrho=\tfrac13$ of uniform distribution, recalling the triangle $C=\{p_0p_1p_2\}$, one can prove that $C\cup G_{\rho,\epsilon}C\subset \overline{A_2\cup B}$ is an IUP, and hence an AsIUP, of $G_{\rho,\epsilon}$ when $\epsilon$ is close enough to $\tfrac12$. Actually, since the point $p_0$ depends on $\epsilon$ when $\varrho\neq \tfrac13$, one needs to adapt the proof of Proposition \ref{PROHd} in this case, using that $\inf_\epsilon \mathrm{dist}(p_0,\overline{A_2}\cap \overline{B})>0$. 

\section{Projection procedure in the case of other permutation groups of $N-1$ coordinates}\label{A-ALTERNAT}
The procedure in Section \ref{S-PROJECT} extends to the case where $F$ commutes with every element of $\Pi_{i_1,\cdots ,i_{N-1}}$, for every $(N-1)$-uple $\{i_1,\cdots ,i_{N-1}\}$. Here, we consider two cases, namely when $\{i_1,\cdots ,i_{N-1}\}=\{2,\cdots ,N\}$ and when $\{i_1,\cdots ,i_{N-1}\}=\{1,\cdots ,N-2,N\}$. We compute the expressions of the corresponding elements, and in particular of the symmetries $\sigma_\Sigma$ and $\phi_N\circ \sigma_\Sigma\circ \phi_N$ associated with the inversion of coordinates sign. Obviously, when $F$ commutes with every permutation in $\Pi_N$, any of these cases can be selected for the reduction procedure. 
 
\subsection{Case of commutation with the permutations of $\{u_i\}_{i=2}^N$}
Here, we assume that $F:\T^N\circlearrowleft$ commutes with every permutation of the coordinates $\{u_i\}_{i=2}^N$. Then, the same statements as in Section \ref{S-PROJECT} hold with the following definitions
\[
D_\ast^N=\left\{u\in [0,1)\times \R^{N-1}\ :\ u_i- u_j\in\R\setminus\Z,\ \forall i\neq j\in [1,N]\ \text{and}\ 0<u_i-u_1<1,\ \forall i\in [2,N]\right\},
\]
\[
(P u)_i=u_i+\lceil u_1-u_i\rceil-\lfloor u_1\rfloor,\quad \forall i\in [1,N],
\]
and 
\[
I_N =\left\{u\in [0,1)\times \R^{N-1}\ :\ u_1<u_2<\cdots <u_{N-1}<u_N<u_1+1\right\}.
\]
The transformation $\pi_u$ is the permutation of the last $N-1$ coordinates that sends $u\in D_\ast^N$ to $I_N$.  

Moreover, the conjugated transformation $\Sigma=P\circ S\circ P^{-1}$ induced by the inversion of coordinate signs $S=-\text{Id}|_{\T^N}$ reads
\[
(\Sigma u)_i=2-\delta_{i,1}-\delta_{u_1,0}-u_i,\ \forall i\in [1,N].
\]
which has proper representation $\sigma_\Sigma$ on $I_N$ given by   
\[
(\sigma_\Sigma u)_i=\left\{\begin{array}{ccl}
1-\delta_{u_1,0}-u_1&\text{if}&i=1\\
2-\delta_{u_1,0}-u_{N-i+2}&\text{if}&i\in [2,N]
\end{array}\right.
\]
which differs from the expression \eqref{INVSYM}. The conjugated transformation $\phi_N\circ \sigma_\Sigma\circ \phi_N$ (see Section \ref{S-COUPLEDM} for the expression of $\phi_N$) writes 
\[
(\phi_N\circ \sigma_\Sigma\circ \phi_N x)_i=\left\{\begin{array}{ccl}
1-(x_1+\cdots+x_{N-1})&\text{if}&i=1\\
x_{N-i+1}&\text{if}&i\in [2,N-1]\\
2-\delta_{(\phi_N^{-1}x)_1,0}-(\phi_N^{-1}x)_{2}&\text{if}&i=N
\end{array}\right.
\]
Notice that $\phi_N\circ \sigma_\Sigma\circ \phi_N$ is a skew-product map whose base map acts on the first $N-1$ coordinates $\{x_i\}_{i=1}^{N-1}$.

\subsection{Case of commutation with the permutations of $\{u_i\}_{i=1}^{N-2}\cup\{u_N\}$}\label{A-ALTERNATSUB2}
Assume that $F:\T^N\circlearrowleft$ commutes with every permutation of the $N-1$ coordinates $\{u_i\}_{i=1}^{N-2}\cup\{u_N\}$. Then, the same statements as in Section \ref{S-PROJECT} hold with the following definitions
\[
D_\ast^N=\left\{u\in \R^{N-2}\times [0,1)\times\R \ :\ u_i- u_j\in\R\setminus\Z,\ \forall i\neq j\in [1,N]\ \text{and}\ 
\left\{\begin{array}{l}
0<u_{N-1}-u_i<1\ \text{if}\ i\in [1,N-2]\\
0<u_N-u_{N-1}<1
\end{array}\right.\right\},
\]
\[
(P u)_i=\left\{\begin{array}{ccl}
u_i+\lfloor u_{N-1}-u_i\rfloor -\lfloor u_{N-1}\rfloor&\text{if}&i\in [1,N-1]\\
u_N+\lceil u_{N-1}-u_N\rceil -\lfloor u_{N-1}\rfloor&\text{if}&i=N
\end{array}\right.
\]
and 
\[
I_N =\left\{u\in \R^{N-2}\times [0,1)\times\R\ :\ u_1<u_2<\cdots <u_{N-1}<u_N<u_1+1\right\}.
\]
In this case, the transformation $\pi_u$ that sends $u\in D_\ast^N$ to $I_N$ is more involved than in the previous cases. It can be described as a two-step process. First, let $\pi^{(1)}$ be the permutation of the first $N-2$ coordinates such that 
\[
(\pi^{(1)} u)_1<(\pi^{(1)} u)_2<\cdots <(\pi^{(1)} u)_{N-2}.
\] 
If we have $u_N<(\pi^{(1)} u)_1+1$, then $\pi^{(1)} u\in I_N$ and $\pi_u=\pi^{(1)}$. In order to describe the other case, assume that $u\in D_\ast^N$ is such that $u_1<\cdots <u_{N-2}$ and $u_1+1<u_N$. Then define
\[
j=\max\{i\in [1,N-2]\ :\ u_i<u_{N}-1\},\quad\text{and}\quad (\pi^{(2)}u)_i=\left\{\begin{array}{ccl}
u_{i+1}&\text{if}&i\in [1,j-1]\\
u_N-1&\text{if}&i=j\\
u_i&\text{if}&i\in [j+1,N-1]\\
u_1+1&\text{if}&i=N
\end{array}\right.
\]
Clearly, we have $\pi^{(2)}u\in I_N$ and hence $\pi_u=\pi^{(2)}$ in this case, so that $\pi_u=\pi^{(2)}\circ \pi^{(1)}$ for an arbitrary $u\in D_\ast^N$ (letting $\pi^{(2)}=\mathrm{Id}$ when $u_1<\cdots <u_{N-2}$ and $u_N<u_1+1$). 
 
Moreover, the conjugated transformation $\Sigma=P\circ S\circ P^{-1}$ induced by the inversion of coordinate signs $S=-\text{Id}|_{\T^N}$ reads
\[
(\Sigma u)_i=\left\{\begin{array}{ccl}
\delta_{i,N-1}-\delta_{u_{N-1},0} -u_i &\text{if}&i\in [1,N-1]\\
2-\delta_{u_{N-1},0}-u_N&\text{if}&i=N
\end{array}\right.
\]
which has proper representation $\sigma_\Sigma$ on $I_N$ given by   
\[
(\sigma_\Sigma u)_i=\left\{\begin{array}{ccl}
-\delta_{u_{N-1},0}-u_{N-2-i}&\text{if}&i\in [1,N-3]\\
1-\delta_{u_{N-1},0}-u_N&\text{if}&i=N-2\\
1-\delta_{u_{N-1},0}-u_{N-1}&\text{if}&i=N-1\\
1-\delta_{u_{N-1},0}-u_{N-2}&\text{if}&i=N
\end{array}\right.
\]
The conjugated transformation $\phi_N\circ \sigma_\Sigma\circ \phi_Ni$ (see Section \ref{S-COUPLEDM} for the expression of $\phi_N$) writes
\[
(\phi_N\circ \sigma_\Sigma\circ \phi_N x)_i=\left\{\begin{array}{ccl}
x_{N-3-i}&\text{if}&i\in [1,N-4]\\
1-(x_1+\cdots +x_{N-1})&\text{if}&i=N-3\\
x_{N-1}&\text{if}&i=N-2\\
x_{N-2}&\text{if}&i=N-1\\
1-\delta_{(\phi_N^{-1}x)_{N-1},0}-(\phi_N^{-1}x)_{N-2}&\text{if}&i=N
\end{array}\right.
\]
which again shows that $\phi_N\circ \sigma_\Sigma\circ \phi_N$ is a skew-product map whose base acts on the first $N-1$ coordinates. Moreover, for $N=3$, this base map is simply the reflection $(x_1,x_2)\mapsto (x_2,x_1)$ with respect to the diagonal.
\end{document}